\documentclass[a4paper, 11pt]{article}

\usepackage{placeins,graphicx,diagbox }
\usepackage{etex}

\usepackage{wrapfig} 
\usepackage{amsmath,tikz,pgfplots,mathtools,float}
\usepackage{booktabs,caption}
\usepackage[flushleft]{threeparttable}
\usetikzlibrary{arrows,shapes,snakes,shadows,positioning,automata,patterns}
\usetikzlibrary{trees,decorations.pathmorphing,decorations.markings}
\usepackage{amsmath}
\usepackage{bm}
\usepackage{subfigure}
\usepackage{multirow}
\usepackage{amssymb,latexsym}  
\usepackage{extarrows}

\usepackage{amsmath,algorithm} 
\usepackage{amssymb}  
\usepackage{color}
\newcounter{algo}

\def\T{^{\rm\tiny T}}

\usepackage{verbatim}
\usepackage{ulem}
\newtheorem{theorem}{Theorem}
\newtheorem{lemma}{Lemma}
\newtheorem{definition}{Definition}

\newtheorem{proposition}{Proposition}

\newtheorem{remark}{Remark}

\title{\bf Data Rates for Network Linear  Equations}

\author{Jinlong Lei, Peng Yi, Guodong Shi, and Brian D. O. Anderson
\thanks{A preliminary version of this work will appear in the Proceedings of the IEEE Conference on Decision and Control \cite{lei2018datatare}. J. Lei is with the Department of Industrial and Manufacturing Engineering, Pennsylvania 	State University, University Park 16802, PA, USA
(email:jxl800@psu.edu); P. Yi is with the Department of Electrical and Systems Engineering, Washington University in St. Louis, 1 Brookings Drive,
St. Louis, MO 63130, USA (email: yipeng@amss.ac.cn); G. Shi  is with the Research School of Engineering, The Australian National University, ACT 0200, Canberra, Australia (email: guodong.shi@anu.edu.au); B. D. O. Anderson is with the Hangzhou Dianzi University, Hangzhou, China,   the Research School of Engineering, Australian National University, Canberra, ACT 0200, Australia (email: brian.anderson@anu.edu.au). The work of Anderson, was supported by the Australian Research Council (ARC) under grant DP-160104500, and by Data61-CSIRO.
}
}
\date{}

\begin{document}

\maketitle

\begin{abstract}
In this paper, we study network linear equations subject to digital communications with a finite data rate, where each node is associated with one equation from a system of linear equations. Each node holds a dynamic state and interacts with its neighbors through an undirected connected graph, where along each link the pair of nodes share information. Due to the data-rate constraint,  each node builds an encoder-decoder pair, with which it produces transmitted message with a zooming-in
 finite-level uniform  quantizer and also generates estimates of its neighbors' states from the received signals.
 We then propose a distributed quantized  algorithm and show that  when the network linear equations admit a unique  solution,
 each node's state is driven to that solution exponentially.
  We further establish the asymptotic rate of  convergence, which shows that a  larger   number of quantization
levels  leads to a faster convergence rate but is still fundamentally bounded by the inherent network structure and the linear equations. When a unique least-squares solution exists, we show that the   algorithm can  compute such a solution with a suitably   selected time-varying step size inherited from the encoder and zooming-in quantizer dynamics. In both cases,  a minimal data  rate is shown to be enough for  guaranteeing the desired convergence when the step sizes are properly chosen. These results assure the applicability  of various network linear equation solvers  in the literature when   peer-to-peer communication is digital.
\end{abstract}


\section{Introduction}
The pursuit of resilient and scalable solutions for the control and optimization of large-scale network systems has been one of the central themes in the field of systems and control in the past decade \cite{tsi,jad03}. For a group of interconnected agents (nodes), sensing and decision making can be carried out individually based on the information flow across the interconnections (links), under which collective goals such as consensus, formation and estimation can be achieved \cite{xiao04,Rabbat2010}. These distributed protocols provide resilience in the sense that nodes and links can join and leave the network without significantly affecting the performance of the network; they also provide   scalability compared to centralized solutions because individual node sensing and decision are often quite simple.
 Simultaneously, control theory has embraced to a much greater degree than previously  graph theory, communication theory, and complexity analysis, leading to many celebrated results for both theories and applications \cite{magnusbook}.

Particularly, systems of linear algebraic equations, as one of  primary computation tasks, can be naturally defined over a network in the way that each node holds one or a few of the linear equations \cite{mou13}. Network linear equations also arise from resource allocation problems when node cost functions are quadratic, see \cite{rabbat2004,nedic09,yi2016initialization}. In the context of parallel computation, computer scientists aimed to develop algorithms that eventually compute part entries of the solutions \cite{margaris2014parallel,saad1999distributed,anderson1997,mehmood05,lei2015distributed}. On the other hand, in view of distributed gradient optimization \cite{nedic10,elia,lei2016primal},  distributed algorithms that compute the entire solution vector at each node were also proposed for both discrete-time  and continuous-time node dynamics \cite{lu12,mou13,Mou-TAC-2015,Shi-TAC-LAE,Yang-IFAC-2017,Yang-CDC-2017,brian15,jadbabaie15,asuman14}. In fact, when exact solutions exist for the linear equations, such first-order  distributed solvers were generalized versions of the so-called alternation projection algorithms pioneered by von Neumann \cite{jvn49,nedic10,shitac}. When no exact solution exists and one considers least-squares solutions, higher-order algorithms or algorithms using properly selected square-summable diminishing step-sizes are needed \cite{Yang-IFAC-2017,Yang-CDC-2017}.

In this paper, we consider network linear equation solvers subject to digital node communications where only   {\it a  finite data rate is available}  \cite{Brockett2000,basar,carli2009,girish2011,hong2016,litao2011}. We use the convenient notion that each node holds one equation from a system of linear equations  with $m$ unknown variables. The nodes aim to reach consensus on the solution of the linear equations. The nodes interconnection is described by an undirected connected graph, where along each link  the neighboring nodes exchange  information constrained by  a limited data rate measured in bits.  Each node builds
an encoder-decoder pair with the help of a zooming-in finite-level uniform quantization function,  and is   equipped with a dynamical internal encoder state co-evolving with the node states. At each  step, each node's encoder produces a quantized message with the node state and the current internal encoder state, which will be transmitted to its neighbors through the digital communication link.  After receiving the quantized information from the neighbors, each node then decodes/estimates its neighbors' states, based on which its own  state is updated with the proposed algorithm. We have established the following results:

\noindent (i)  When the network linear equation admits a unique exact solution, we show that the proposed encoder-decoder powered algorithm drives each node state to that solution asymptotically with an exponential convergence rate based on   merely $m$   bits information exchange between each pair of adjacent agents.
Furthermore, we give an explicit form of the asymptotic rate of  convergence, which  is
related to  the scale and the synchronizability  of the network, the number of quantization levels, the dimension of the unknown variable,   and   the observation matrix.  It is shown that   a higher convergence rate is possible with higher  data rates but is fundamentally bounded by the inherent network structure and the linear equations.

\noindent  (ii)  When the network linear equation admits a unique least-squares solution, we show that the same encoder-decoder pair enables the algorithm to compute such a solution with a time-varying step-size that comes from the dynamics of encoder internal states. Again,  a  data  rate of $m$ bits per step can deliver such a convergence result, and an explicit form of the asymptotic convergence rate  is established.

These results serve as assurance of the practical use of the various network linear equation solvers  when digital point-to-point communications are subject to round-up errors. Generalizations to the scenarios where the solutions of the linear equations are not unique for both exact and least-squares cases are possibly along the same line of analysis, but are not included in the current paper for the ease of presentation. We also note that our results are closely related to the work on distributed optimization algorithms with quantized communication \cite{rabbat2005,peng2014}. However, new challenges for network linear equations arise compared to distributed optimization framework, although the problem appears to be a   special case of quadratic  program at first glance, lie in that  gradients of the quadratic function associated with each node cannot be assumed to be globally bounded  a priori, a key technical assumption for the convergence results of distributed (sub)gradient optimization \cite{nedic10,peng2014}.

A preliminary version of the results  will  be presented at  the  IEEE CDC  in 2018 \cite{lei2018datatare}.
Current manuscript compared to  \cite{lei2018datatare} makes the following improvements and extensions:
 (i)  we future specify how the rate of convergence is influenced by the  quantization levels, the scale and the synchronizability of the network, the  variable dimension $m$ as well  as the problem structure;
 (ii) we carry out more simulations to discuss how data rate influences algorithm parameter selection, and thereby, influences the converge rate; (iii) we also compare convergence rates for different types of communication graphs,  and
 give the  completed proofs of all results.
The remainder of this paper is organized as follows. Section \ref{sec:problem} defines the network linear equation, introduces the node encoders and decoders, and develops a distributed quantized  algorithm.     Section \ref{sec:exact} presents the   exact solver along with its convergence analysis and numerical examples. Section \ref{sec:leastsquares} further investigates the least-squares case. Finally, concluding remarks are given in Section \ref{sec:conclusions}.

\medskip

\noindent{\em Notation and Terminology}. All vectors are column vectors and denoted by bold, lower case letters, i.e., $\mathbf{a},\mathbf{b},\mathbf{c}$,  etc.; matrices are denoted with bold, upper case letters, i.e.,  $\mathbf{A},\mathbf{B},\mathbf{C}$,  etc.;  sets are denoted with $\mathcal{A},\mathcal{B},\mathcal{C}$, etc. Depending on the argument, $|\cdot|$ stands for the absolute value of a real number or the cardinality of a set.
 The Euclidean norm of a vector is denoted as $\|\cdot\|$. $\otimes$ denotes the Kronecker  product. An undirected graph  is an ordered pair of two sets denoted by  $\mathcal{G}=\{ \mathcal{V},\mathcal{E}\}$ where $\mathcal{V}=\{1,\dots,N\}$  is a finite set of  vertices (nodes), and each element in $\mathcal{E}$ is an unordered  pair of two distinct  nodes in $\mathcal {V}$, called an edge.  A    path in $\mathcal{G}$ with length $p$ from $v_1$ to $v_{k+1}$ is a  sequence of distinct nodes, $v_1v_2\dots v_{p+1}$, such that  $(v_m, v_{m+1}) \in \mathcal{E}$, for all $m=1,\dots,p$. The graph $\mathcal{G}$ is termed {\it   connected} if for any two distinct nodes $i,j\in\mathcal{V}$, there is a  path between them. The neighbor set of node $i$, denoted $\mathcal{N}_i$, is defined as  $\mathcal{N}_i=\{ j\in \mathcal{V}: (i,j)\in \mathcal{E}\}$.
 Define the degree matrix $ \mathbf{  D_e}=  {\rm diag} \{ | \mathcal{N}_1 | ,\dots, |\mathcal{N}_N |\} $
 and the adjacency matrix $\mathbf{ A}$, where $[\mathbf{ A}]_{ij}=1$  if $j\in \mathcal{N}_j$ and
 $[\mathbf{ A}]_{ij}=0$  otherwise.  Then $\mathbf{ L}=\mathbf{ D_e}-\mathbf{ A}$ is the Laplacian  matrix of the graph $\mathcal{G}$.

\section{Problem Statement and Algorithm Design}\label{sec:problem}


\subsection{Linear Equations over Networks}
Consider the following  linear algebraic equation:
\begin{align}\label{LinearEquation}
\mathbf{z}=\mathbf{H} \mathbf{y}
\end{align}
with respect to unknown variable $\mathbf{y}\in\mathbb{R}^m$, where $\mathbf{H}\in\mathbb{R}^{N\times m}$ and $\mathbf{z}\in\mathbb{R}^N$. The equation (\ref{LinearEquation}) has a unique exact solution if ${\rm rank}(\mathbf{H})=m$ and $\mathbf{z}\in {\rm span}(\mathbf{H})$; an infinite set of  solutions if ${\rm rank}(\mathbf{H})<m$ and $\mathbf{z}\in {\rm span}(\mathbf{H})$; and no exact solutions   if  $\mathbf{z}\notin {\rm span}(\mathbf{H})$.  When no exact solution exists, a   least-squares   solution of  \eqref{LinearEquation}  can be defined via the following  optimization problem:
\begin{align}\label{LS}
\min_{\mathbf{y}\in\mathbb{R}^m} \big\|\mathbf{z}-\mathbf{H}\mathbf{y}\big\|^2,
\end{align}
which yields a unique solution $\mathbf{y}^\star=(\mathbf{H}\T\mathbf{H})^{-1}\mathbf{H}^T\mathbf{z}$ if ${\rm rank}(\mathbf{H})=m$.

We denote by
$$
\mathbf{H}= \begin{pmatrix}
  \mathbf{h}_1\T \\
 \mathbf{h}_2\T \\
  \vdots  \\
\mathbf{h}_N\T
 \end{pmatrix}, \quad \mathbf{z}= \begin{pmatrix}
  z_1 \\
 z_2 \\
  \vdots  \\
z_N
 \end{pmatrix},
$$
where $\mathbf{h}_i\in \mathbb{R}^m$ with  $\mathbf{h}_i\T$ being the $i$-th row vector of $\mathbf{H}$.

Consider a network with  $N$ nodes  indexed as $\mathcal{V}=\big\{1,\dots,N\big\}$, where node $i$ has access to the value of $\mathbf{h}_i$ and ${z}_i$ without the knowledge of $\mathbf{h}_j$ or ${z}_j$ from other nodes. The nodes interaction   is described by a connected undirected  graph $\mathcal{G}=\{ \mathcal{V},\mathcal{E}\}$  with the corresponding Laplacian matrix denoted by $\mathbf{L}$.  Time is slotted at $k=0,1,2,\dots$. Node $i$  at time $k$ holds an estimate $\mathbf{x}_i(k) \in\mathbb{R}^m$
for the solution to equation  \eqref{LinearEquation} and exchanges information with its neighbors.

As the Euler approximation of the so-called ``consensus + projection" flow proposed in \cite{Shi-TAC-LAE}, the following algorithm is  an efficient distributed linear equations solver with a  discrete recursion.
\begin{equation}\label{1}
\begin{split}
\mathbf{x}_i(k+1) =\mathbf{x}_i(k)&+h \Big[   \sum_{j\in \mathcal{N}_i}\big(\mathbf{x}_j(k)-\mathbf{x}_i(k)\big)    -  \gamma(k)\left(  \mathbf{h}_i\mathbf{h}_i^\top  \mathbf{x}_i(k)- z_i\mathbf{h}_i \right) \Big] .
\end{split}
\end{equation}
It can be easily concluded  from the analysis in \cite{Shi-TAC-LAE,Yang-CDC-2017} that the following statements hold for the algorithm (\ref{1}):
\begin{itemize}
\item When the linear equation (\ref{LinearEquation}) admits an exact solution $\mathbf{y}^*$, it drives each $\mathbf{x}_i(k)$ to $\mathbf{y}^*$ exponentially with   $\gamma(k)\equiv \gamma>0$  provided that  $h, \gamma$ are properly chosen.

     \item When the linear equation \eqref{LinearEquation} has no exact solutions, it drives each $\mathbf{x}_i(k)$ to a   least-squares   solution  to \eqref{LS} for small $h$ and   $\gamma(k)=1/k$.
\end{itemize}
\medskip

It is clear that in the algorithm (\ref{1}), nodes need to exchange  their exact state values for the execution of the update. The aim of this paper is to develop algorithms that overcome such a constraint using {\em quantized} node communications,
and to explore  the corresponding   convergence properties  with minimal data rate statements.


\subsection{ Distributed Quantized  Algorithm}

Suppose that the communication channels  corresponding to each edge in the network have a limited  capacity or  a finite bandwidth.  As such, real-valued data should be quantized before transmitting.
Thus, we propose  a distributed quantized algorithm, in which  each node is associated with an encoder while  its neighbors possess a corresponding decoder.
Let us begin by introducing   a  uniform quantization function  $Q_K(\cdot)$.

\begin{definition}[Quantization Function]  A  standard uniform quantizer is given by the function   $ Q_K(\cdot):\mathbb{R}\to \{-K,\dots,-1,0,1,\dots,K\}$  where
 \begin{align}\label{quantizer}
 Q_K(z)=\begin{cases}  & 0,~ {\rm if~}  -1/2\leq z \leq 1/2, \\
& i,~{\rm if~}  {2i-1\over 2}<z\leq {2i+1\over 2}, ~i=1,\dots, K ,\\
& K,~ {\rm if~}  z > {2K+1 \over 2}, \\
& - Q_K(-z) ,\quad {\rm if~}   z>-1/2.
\end{cases}
\end{align}
\end{definition}
There is no need to send any information if the output of the quantizer  is zero, so,  for a $2K+1$-level
quantizer,  the communication channel  is required to be capable of transmitting $  \lceil \log_2(2K ) \rceil$ bits.
 With slight abuse of notation, we  define $ Q_K(\mathbf{a})$
for a vector $\mathbf{a}=(a_1,\dots, a_m)^T \in \mathbb{R}^m$ by
$$ Q_K(\mathbf{a})=(  Q_K(a_1),\dots,  Q_K(a_m))^T.$$

Next, we propose  an   encoder-decoder pair   for each node to  quantize its state,   and  to estimate the neighbors' states.
Suppose the nodes have a global scaling function $s(k)$.  We still use $\mathbf{x}_i(k)$  to denote the un-quantized state of   node  $i$  at time $k$, whose update will be specified at a later stage.

 \begin{algorithm} [H]
\caption*{\bf Encoder}
Node  $j\in \mathcal{V}$  recursively  generates  $m$-vector    quantized  outputs  $\{\mathbf{q}_j(k)\}$
 and  $m$-vector   internal states $\{\mathbf{b}_j(k)\}$ from the exact  $m$-vector   state sequence  $\{\mathbf{x}_i(k)\}$   as follows for any $k\geq 1$:
\begin{equation} \label{encoder}
 \begin{split}
&\mathbf{q}_j(k)\triangleq  Q_K\left({1\over s(k-1)}(\mathbf{x}_j(k)-\mathbf{b}_j(k-1))\right),\\
&\mathbf{b}_j(k)\triangleq  s(k-1)\mathbf{q}_j(k)+\mathbf{b}_j(k-1),
\end{split}
\end{equation}
where the initial value $\mathbf{b}_j(0)=0$.
\end{algorithm}

%
%

\begin{remark} Note that $\mathbf{b}_j(k)$ is a one-step predictor, and the encoder is a difference encoder with
a zooming-in scaling $s(k)$ that quantizes  the  prediction error  $ \mathbf{x}_j(k)-\mathbf{b}_j(k-1) $  rather than the state $ \mathbf{x}_j(k) $. Generally speaking, the amplitude of the prediction error is smaller than that of the state itself, so it can be represented by
fewer bits.  We use  the scaling function  $s(k)$  to zoom-in  each node's  prediction error and require that $s(k)$  decay gradually
to make the quantizer persistently excited, such that the nodes gradually increase the accuracy of state recovery  of their neighbors. On the other hand, $s(k)$ should be large enough such that the quantizer will not be saturated, in which case the quantization error is bounded. We revisit subsequently the issue of avoidance of saturation.
\end{remark}

Node $j \in \mathcal{V}$ at time $k$ sends its  quantized  output    $\mathbf{q}_j(k)$   to  its neighboring nodes  $i\in \mathcal{N}_j$, which then  recovers node $j$'s state  using the  decoder defined as follows.

 \begin{algorithm} [H]
\caption*{\bf Decoder}
 When node $i\in \mathcal{N}_j$  receives    the quantized data $\mathbf{q}_j(k)$  from node $j$, a decoder recursively generates    an estimate $ \hat{\mathbf{x}}_{ij}(k)$  for  $\mathbf{x}_j(k)$ by the following  for any $k\geq 1$:
\begin{equation}\label{decoder}
  \hat{\mathbf{x}}_{ij}(k) \triangleq s(k-1)\mathbf{q}_j(k)+\hat{\mathbf{x}}_{ij}(k-1),
  \end{equation}
where the initial value   $\hat{\mathbf{x}}_{ij}(0) \triangleq {0}$.
\end{algorithm}

Based on the encoder-decoder pair defined in  \eqref{encoder} and \eqref{decoder}, motivated by \eqref{1}, we now propose the following distributed linear equation solver with quantized node communication.

 \begin{algorithm} [H]
\caption{Distributed quantized algorithm}  \label{alg1}
\begin{eqnarray}\label{quantized-do}
\mathbf{x}_i(k+1) = &\mathbf{x}_i(k)+h    \Big[\sum_{j\in \mathcal{N}_i}  \big(\hat{\mathbf{x}}_{ij}(k)-\mathbf{b}_i(k)\big)   -\gamma(k) \left(  \mathbf{h}_i\mathbf{h}_i^\top \mathbf{x}_i(k)- z_i\mathbf{h}_i \right) \Big] .
\end{eqnarray}
\end{algorithm}

%
%

It is worth noting that   the  difference between  \eqref{1} and \eqref{quantized-do} lies in  the fact  that
  the exact state $ \mathbf{x}_{j }(k)$ is  used in \eqref{1}
while    $\hat{\mathbf{x}}_{ij}(k)$ is used in  \eqref{quantized-do}.
It is clear that Algorithm  \ref{alg1} merely relies on quantized node communication   since $\mathbf{q}_j(k)$ takes values in the alphabet $\{-K,\dots,-1,0,1,\dots,K\}$ only.
 From the second equation of \eqref{encoder}, using Equ. \eqref{decoder} and    the assumed initial conditions of zero for $\hat{\mathbf{x}}_{ij}(0)$ and  $\mathbf{b}_j(0)$,   we   have the following for any $k\geq 0$:
\begin{equation}\label{equiv}
\hat{\mathbf{x}}_{ij}(k)=\mathbf{b}_j(k), \quad  \forall j \in \mathcal{V},\  \forall  i \in \mathcal{N}_j.
\end{equation}

\section{Exact Solutions}\label{sec:exact}
In this section, we  consider  Algorithm \ref{alg1} and investigate  the case that   equation   \eqref{LinearEquation} has a unique solution. We establish the convergence results regarding the  quantization levels  along with
the rate analysis, and demonstrate the results with numerical simulations.

\subsection{Convergence Result}

We impose the following  assumptions.

\noindent{\bf A1}  There exists a unique solution $\mathbf{y}^*$, i.e.,  ${\rm rank}(\mathbf{H})=m$ and $\mathbf{z}\in {\rm span}(\mathbf{H})$.

\medskip

\noindent{\bf A2}  $\max_i \|\mathbf{x}_i(0)\|_{\infty}\leq C_x$ and  $\max_i \|\mathbf{x}_i(0)-\mathbf{y}^* \|_{\infty}\leq C_w$ for some positive constants  $C_x$ and $C_w$.

\medskip

\noindent{\bf A3}  $\gamma(k)\equiv 1  $, and  $s(k)\triangleq s(0)  \alpha^k~\forall k\ge 0$ for some $s(0)>0 $ and $\alpha \in (0,1).$

We now introduce a few useful notations as follows:
 \begin{equation}\label{def}
\begin{split}
&   \mathbf{H_{d}}  \triangleq  {\rm diag}  \left  \{ \mathbf{h}_1\mathbf{h}_1^\top ,\dots,  \mathbf{h}_N\mathbf{h}_N^\top  \right \} \in \mathbb{R}^{mN\times mN},\\&   \mathbf{F_d}\triangleq   \mathbf{L}\otimes \mathbf{I}_m + \mathbf{H_d},
~   \rho_h \triangleq  1-h  \lambda_{\min}(\mathbf{F_d}),
    \end{split}
\end{equation}
 where  $ \lambda_{\min}(\mathbf{F_d})$   denotes  the smallest   eigenvalue  of $  \mathbf{F_d} $.
 Note that both the Laplacian matrix $\mathbf{L}$   and  the matrix $  \mathbf{H_d}$ are positive semidefinite. With the assumption {\bf A1} and the condition that   the  undirected  graph $\mathcal{G}$ is connected, the matrix $\mathbf{F_d}$
 turns  out to be positive definite \cite[Lemma 9]{Shi-TAC-LAE}, and hence all eigenvalues of  $\mathbf{F_d}$ is positive.
  The  eigenvalues of  $\mathbf{L}$  in an ascending order
are denoted by  $0=\lambda_1(\mathbf{L})<\lambda_2(\mathbf{L})\leq \dots \leq \lambda_N(\mathbf{L}) .$ Let  $h \in \left(0, {2\over \lambda_{\min}(\mathbf{F_d})+\lambda_{\max}(\mathbf{F_d})}\right)$ and   $\alpha \in (  1-h  \lambda_{\min}(\mathbf{F_d}),1)$, and set
 \begin{equation}\label{def-M}
\begin{split}
& M(\alpha, h)\triangleq{ 1+2hd^* \over 2 \alpha}+
{h^2 \sqrt{mN} \lambda_N(\mathbf{L}) \lambda_{\max}(\mathbf{F_d})\over 2 \alpha(\alpha-\rho_h)} ,  {\rm~and ~} \mathcal{K}(\alpha, h )\triangleq \Big \lceil M(\alpha, h) -{1\over 2}\Big \rceil  ,
  \end{split}
\end{equation}
where $d^*=\max_{i} | \mathcal{N}_i| $ denotes  the degree of $\mathcal{G},$  and  $ \lambda_{\max}(\mathbf{F_d})$   denotes  the largest   eigenvalue  of $  \mathbf{F_d} $.

 We  now  begin to investigate the convergence properties of Algorithm  \ref{alg1}  as an exact solver
for the network linear equation \eqref{LinearEquation}.

 \begin{proposition}[Non-Saturation] \label{prp1}
Let  {\bf A1}, {\bf A2} and {\bf A3}  hold.
Consider Algorithm  \ref{alg1}, where
\begin{align*} h \in \left(0, {2\over \lambda_{\min}(\mathbf{F_d})+\lambda_{\max}(\mathbf{F_d})}\right) {~\rm and ~}  \alpha \in (1-h \lambda_{\min}(\mathbf{F_d}),1) .\end{align*}  Then for any   $K \geq \mathcal{K}(\alpha, h )$, the quantizer will never be saturated
provided that $s(0)$ satisfies   
 \begin{equation} \label{def-s0}
\begin{split}
s(0)>   \max  \Bigg\{  &{ C_x+ h  \|   \mathbf{H_d}  \|_{\infty} C_w \over K+{1\over 2}},  {2   (\alpha-\rho_h)
\left(  \rho_h C_w  + hC_x  \lambda_N(\mathbf{L}) \right) \over h  \lambda_N(\mathbf{L}) } \Bigg \} .
\end{split}
\end{equation}  
\end{proposition}
 Proposition  \ref{prp1} with the proof  deferred to Section \ref{proof-thm01} establishes the nonsaturation of the uniform quantizer, based on which  the  following theorem  with the proof   given in Section \ref{proof-thm1} shows  the asymptotic convergence of the generated sequences to the unique exact solution.

\begin{theorem} [High Data Rate]\label{thm1}
Suppose {\bf A1}, {\bf A2} and {\bf A3}  hold.
With $\mathbf{F_d}$ as defined in \eqref{def}, let    $h \in \left(0, {2\over \lambda_{\min}(\mathbf{F_d})+\lambda_{\max}(\mathbf{F_d})}\right)$ and  $\alpha \in (1-h \lambda_{\min}(\mathbf{F_d}),1)$.
 Then for any   $K \geq\mathcal{K}(\alpha, h )$, see \eqref{def-M}, along Algorithm  \ref{alg1} there holds
\begin{align}
& \lim_{k \to \infty} x_i(k)=\mathbf{y}^*\quad \forall i \in \mathcal{V} \label{convergence}
\end{align}   provided    $s(0) $  satisfying    \eqref{def-s0}.
The convergence is in fact  exponential with
\begin{align}
\limsup_{k \to \infty }{  \|\mathbf{x}(k)- \mathbf{1}_N \otimes \mathbf{y}^* \|_2 \over \alpha^k}  \leq  {hs(0) \sqrt{mN}  \lambda_N(\mathbf{L}) \over 2  \alpha  (\alpha-\rho_h)}   \label{rate},
\end{align}
where  $\mathbf{x}(k)= col\{ \mathbf{x}_1(k),\dots, \mathbf{x}_N(k)\}
\triangleq(\mathbf{x}_1(k)^T,\dots, \mathbf{x}_N(k)^T)^T$.
\end{theorem}
\begin{remark}
Theorem \ref{thm1} shows that by using a scaling function decaying exponentially and a  uniform quantizer,  Algorithm  \ref{alg1}   can ensure  asymptotic convergence to the unique solution. It is worth pointing out that for any given $\alpha,h$, the  obtained quantization level $\mathcal{K}(\alpha,h)$
is  conservative, while  \eqref{def-M}  gives us some intuition   on the relationship between the number of bits required
and the control gains and the scaling factor. In addition, Theorem \ref{thm1}  gives an estimate of the rate of convergence: the smaller   the scaling factor $\alpha$, the faster the convergence rate  from \eqref{rate}  but  more bits have to be communicated by  \eqref{def-M},  and, if $\alpha \to \rho_h$,   the required number of bits goes to infinity. Thus, an appropriate selection of $\alpha$  amounts to a tradeoff  between the rate of convergence  and the communication overhead.
\end{remark}

From  \eqref{def-M} we   know that for fixed $\alpha$,  the quantization level  $\mathcal{K}(\alpha ,h) $  will tend to infinity as $N\to \infty$.   Since   in practical applications,  the communication channel   usually has  finite bandwidth. To satisfy this requirement, we can use  a fixed number of quantization levels at the cost of slower convergence.
We present the result in  the following  theorem, for which the proof is given in Section \ref{proof-thm2}.

\begin{theorem}[Low Data Rate] \label{thm2}
 Suppose {\bf A1}, {\bf A2}, and {\bf A3}  hold, with $\mathbf{F_d}$ and $M(\alpha,h)$ as defined in \eqref{def} and \eqref{def-M}.
Then  the following hold.\\
(i) For any $K\geq 1,$ $\Xi_K$ is nonempty with \begin{equation}\label{def_omega}
\begin{split}
\Xi_K \triangleq & \Big \{ (\alpha,h ):    h \in \Big(0, {2\over \lambda_{\min}(\mathbf{F_d})+\lambda_{\max}(\mathbf{F_d})}\Big) ,    \alpha\in (1-h  \lambda_{\min}(\mathbf{F_d}) ,1), M(\alpha, h) <  K+{1\over 2}\Big\}.
\end{split}
\end{equation}
(ii)  For any  $K\geq 1$, let $ (\alpha,h )\in \Xi_K$ and  $s(0) $  satisfy \eqref{def-s0}. Then  along  Algorithm  \ref{alg1}  there holds
$ \lim_{k \to \infty} x_i(k)=\mathbf{y}^*~ \forall i \in \mathcal{V}$ at an exponential rate characterized by
$$\limsup_{k\to \infty }{  \|\mathbf{x}(k)- \mathbf{1}_N \otimes \mathbf{y}^* \|_2 \over \alpha^k}  \leq  {hs(0) \sqrt{mN}  \lambda_N(\mathbf{L}) \over 2  \alpha  (\alpha-\rho_h)} .$$
\end{theorem}

\begin{remark}
 From  Theorem  \ref{thm2} it is clear  that we can always design a distributed network linear equation solver  to ensure exponential convergence to exact solution with $3-$levels quantizer (namely, $K=1$), under which each node sends merely  $m$ bits of information ({\bf minimum  number of   bits}) to its  neighbors at each   step.
\end{remark}

 From definition \eqref{def_omega} it is seen that the  set  $\Xi_K$ is defined by three nonlinear
inequalities, for which   an explicit solution of these inequalities  might be difficult to obtain.
Then in   the following proposition with the proof given in Section \ref{proof-lem}, we give an explicit  method  for   choosing  parameters  $(\alpha,h )$
from $\Xi_K$ for any given $K \geq 1$ by  introducing a free parameter $\epsilon \in (0,1)$.
\begin{proposition} \label{lem-rate}
 For any given $K\geq 1 $  and   $\epsilon\in (0,1)$,  define  $ \Xi_{K,\epsilon} \triangleq   \big \{ (\alpha,h ):  \alpha=1-(1-\epsilon) h \lambda_{\min}(\mathbf{F_d}) , h\in (0,h_{K,\epsilon}^*) \big \}  ,  $
where  $
h_{K,\epsilon}^*\triangleq  \min \left\{ {2\over \lambda_{\min}(\mathbf{F_d})+\lambda_{\max}(\mathbf{F_d})},\hat{h}_{K,\epsilon} \right\}  $  with   
\begin{equation}\label{def-hath}
\begin{split} &\hat{h}_{K,\epsilon} \triangleq  2K  \epsilon  \lambda_{\min}(\mathbf{F_d})  \Big(  \sqrt{mN} \lambda_N(\mathbf{L}) \lambda_{\max}(\mathbf{F_d}) +2 \epsilon  \lambda_{\min}(\mathbf{F_d})  d^* + \epsilon  (1-\epsilon) (2K+1)    \lambda^2_{\min}(\mathbf{F_d})    \Big)^{-1}.
\end{split}
\end{equation}
Then we have that
$\Xi_K=\bigcup_{\epsilon \in (0,1)} \Xi_{K,\epsilon} .$
\end{proposition}

 We note from Theorem \ref{thm1}   that the proposed distributed protocol ensures exponential convergence  with parameter $\alpha,$ which   is coupled with another algorithm parameter $h$ while without explicit dependence on the linear equations and the network.
In  the following, we   investigate the asymptotic property  of  $\alpha$ as $N\to \infty$, and give a very compendious expression for   the asymptotic value  of  $\alpha$.
The  proof can be found    in Section \ref{proof-thm-rate}.
\begin{theorem}[Network Scalability] \label{thm-rate}
Adopt the same hypothesis as Theorem \ref{thm2}. Let $K\geq 1$ and $ (\alpha,h )\in \Xi_K$. Then
\begin{align} \label{bd-alpha}
\lim_{N \to \infty} { \inf_{ (\alpha,h)  \in \Xi_K}\alpha \over  {\rm exp} \left( - { K  \lambda^2_{\min}(\mathbf{F_d}) \over 2\sqrt{mN} \lambda_N(\mathbf{L}) \lambda_{\max}(\mathbf{F_d})}  \right)   } =1.
\end{align}
\end{theorem}

\begin{remark}
Theorem \ref{thm-rate}  together with Equ. \eqref{rate}  suggests that  for large network, the highest possible rate of convergence tends to scale according to $  \mathcal{O} \left( {\rm exp} \left( - k \Theta_N K   \right) \right),$
where $$\Theta_N = {{\lambda^2_{\min}(\mathbf{F_d})} \over {2\sqrt{mN} \lambda_N(\mathbf{L}) \lambda_{\max}(\mathbf{F_d})}}$$ is some constant relying  only on the number of nodes, the network structure and the equations.
\end{remark}

\subsection{Numerical Examples}\label{subsec:exact:ls:quantizaiton}

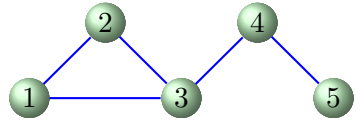
\begin{figure}[!htb]
\begin{center}
\begin{tikzpicture}[->,>=stealth',shorten >=0.3pt,auto,node distance=2cm,thick,
  rect node/.style={rectangle, ball color={rgb:red,0;green,0.2;yellow,1},font=\sffamily,inner sep=1pt,outer sep=0pt,minimum size=14pt},
  wave/.style={decorate,decoration={snake,post length=0.1mm,amplitude=0.5mm,segment length=3mm},thick},
  main node/.style={shape=circle, ball color=green!20,text=black,inner sep=1pt,outer sep=0pt,minimum size=15pt},scale=1]
  \foreach \place/\i in {{
    (-2,0)/1},
   {(-1,1)/2},
   {(0,0)/3},
   {(1,1)/4},
   {(2,0)/5}}
    \node[main node] (a\i) at \place {};

      \node at (-2,0){\rm \color{black}{$1$}};
      \node at (-1,1){\rm \color{black}{$2$}};
      \node at (0,0){\rm \color{black}{$3$}};
      \node at (1,1){\rm \color{black}{$4$}};
      \node at (2,0){\rm \color{black}{$5$}};

  \path[-,blue,thick]               (a1) edge (a2);
  \path[-,blue,thick]               (a1) edge (a3);
  \path[-,blue,thick]               (a2) edge (a3);
  \path[-,blue,thick]               (a3) edge (a4);
  \path[-,blue,thick]               (a4) edge (a5);
 \end{tikzpicture}
\end{center}
\caption{Communication graph.}\label{fig_communication_graph}
\end{figure}

\noindent {\bf Example 1.}
Let the linear equation (\ref{LinearEquation}) be given by
\begin{equation}
\begin{split}\label{example_H_z}
&\mathbf{H}=\left(
    \begin{array}{cc}
      0.5	& -0.1\\
    -0.4  &	0.2\\
   0.3	&  -0.7\\
   0.6	&  0.3\\
   -0.3   &	0.5\\
    \end{array}
  \right),
  \mathbf{z}=\left(
    \begin{array}{c}
    0.2\\
    0.2\\
   -1.8\\
    1.5\\
    1.2\\
    \end{array}
  \right)\end{split}
\end{equation}
which  yields a unique exact solution
\begin{equation*}
 \mathbf{y}^*=
  \left(
    \begin{array}{c}
      1\\
3\\
    \end{array}
  \right).
\end{equation*}
The network structure is  shown in Figure \ref{fig_communication_graph}.

\begin{figure}[htbp]
\centering
\begin{minipage}[t]{0.48\textwidth}
\centering
 \includegraphics[width=3.6in]{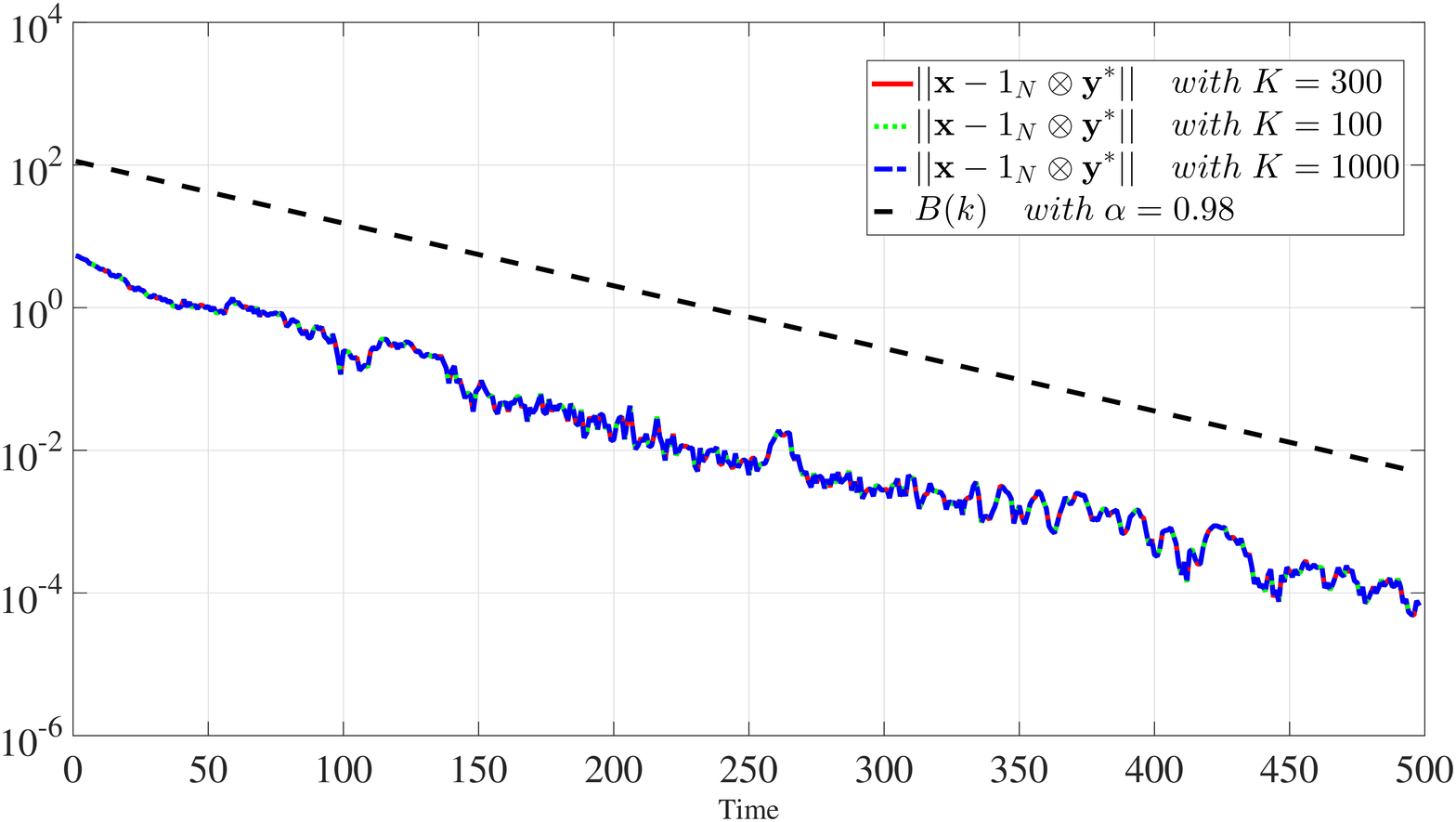}\\
  \caption{\small Trajectories of  $ ||\mathbf{x}(k)-\mathbf{1}_N \otimes \mathbf{y}^* ||_2$ along with the upper bound $B_k$ under $K= 100,300, 1000$.}\label{fig_exact_thm1 1}
\end{minipage}
\begin{minipage}[t]{0.04\textwidth}
\end{minipage}
\begin{minipage}[t]{0.48\textwidth}
\centering
   \includegraphics[width=3.6in]{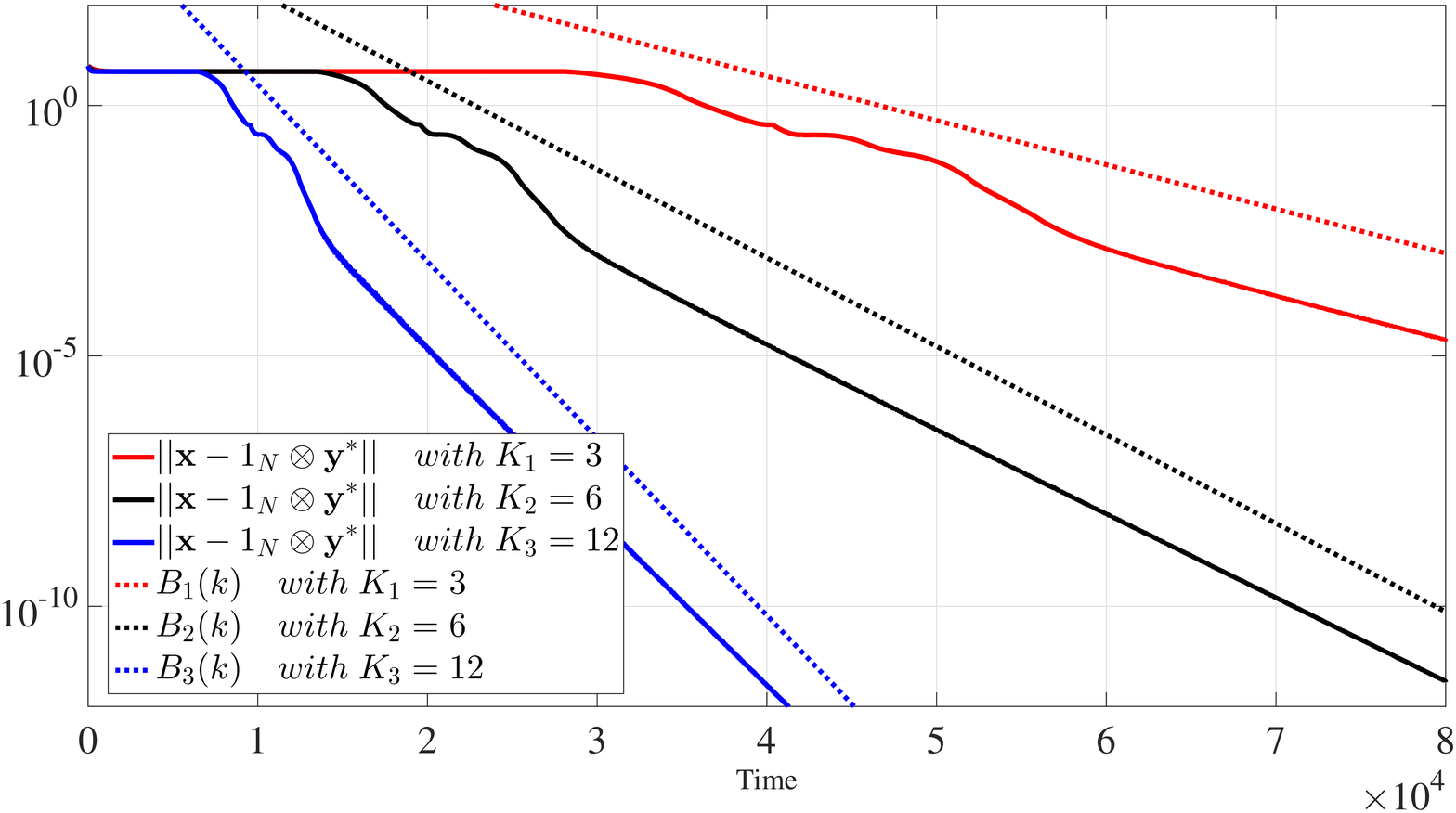}\\
  \caption{\small Trajectories of   $ ||\mathbf{x}(k)-\mathbf{1}_N \otimes \mathbf{y}^* ||_2$ and  $B(k)$ with $K_1=3$, $K_2=6$ and $K_3=12$, respectively.}\label{fig_exact_thm2_1}
\end{minipage}
\end{figure}

\noindent[{\bf Validation of  Theorem \ref{thm1}.}]
Let $h= \frac{1.98}{\lambda_{\min}(\mathbf{F_d})+\lambda_{\max}(\mathbf{F_d})}=0.4215$. Here one can
compute $\rho_h=0.9554$.
Set $\alpha=0.98$  so that  $\mathcal{K}(\alpha,h)=225$.
Let $K$ be $100$, $300$,  $1000$, respectively. We set $s(0)=1$ and implement  Algorithm \ref{alg1}.
Figure \ref{fig_exact_thm1 1} displays the trajectories of $ ||\mathbf{x}(k)-\mathbf{1}_N \otimes \mathbf{y}^* ||_2$   along with  the theoretical  upper bound $B(k)= {hs(0) \alpha^k \sqrt{mN}  \lambda_N(\mathbf{L}) \over 2  \alpha  (\alpha-\rho_h)}   $ given by \eqref{rate}. The trajectory with $K=300$ verifies that Theorem 1 provides a sufficient condition on the data rate to ensure convergence, while the trajectories for $K=100$ and $K=1000$ coincide with that   of $K=300$.  Therefore, it  implies that  (i) with the same algorithm parameters $h,\alpha $, a higher data rate ($K=1000$) cannot guarantee  a faster convergence rate; (ii) there is some degree of conservativeness in the sufficient condition of Theorem \ref{thm1}.

 \noindent[{\bf Validation of   Theorem \ref{thm2}.}] Let $K$  be $K_1=3$, $K_2=6$ and $K_3=12$, respectively. We choose $(\alpha,h)\in \Xi_K$ with Proposition \ref{lem-rate}. Set $\epsilon=0.5$, and we  then  choose
 $  (\alpha_1,h_1)=(0.9998,0.0038) \in \Xi_{K_1,0.5},  $ $(\alpha_1,h_1)=(0.9996,0.0077) \in \Xi_{K_2,0.5}, $ and $(\alpha_3,h_3)=(0.9992,0.0154) \in \Xi_{K_3,0.5}$.
We set $s_1(0)=1500,s_2(0)=1200,s_3(0)=1000$ for $K_1$, $K_2$, $K_3$, respectively, to ensure \eqref{def-s0}.
The trajectories of $ ||\mathbf{x}(k)-\mathbf{1}_N \otimes \mathbf{y}^* ||_2$ under the three sets of parameters are shown in  Figure \ref{fig_exact_thm2_1}, which demonstrates    the convergence of Algorithm \ref{alg1} to the exact solution. A higher data rate allows us to choose  a larger $h$ and a smaller $\alpha$, and therefore, leads to
 a faster convergence rate.   Figure \ref{fig_exact_thm2_1} is also consistent with the  upper bound  of convergence rate $B(k)= {hs(0) \alpha^k \sqrt{mN}  \lambda_N(\mathbf{L}) \over 2  \alpha  (\alpha-\rho_h)}   $ given by \eqref{rate}
 in all three parameter settings.
%

\begin{figure}[htbp]
\centering
\begin{minipage}[t]{0.46\textwidth}
\centering
  \includegraphics[width=3.6in]{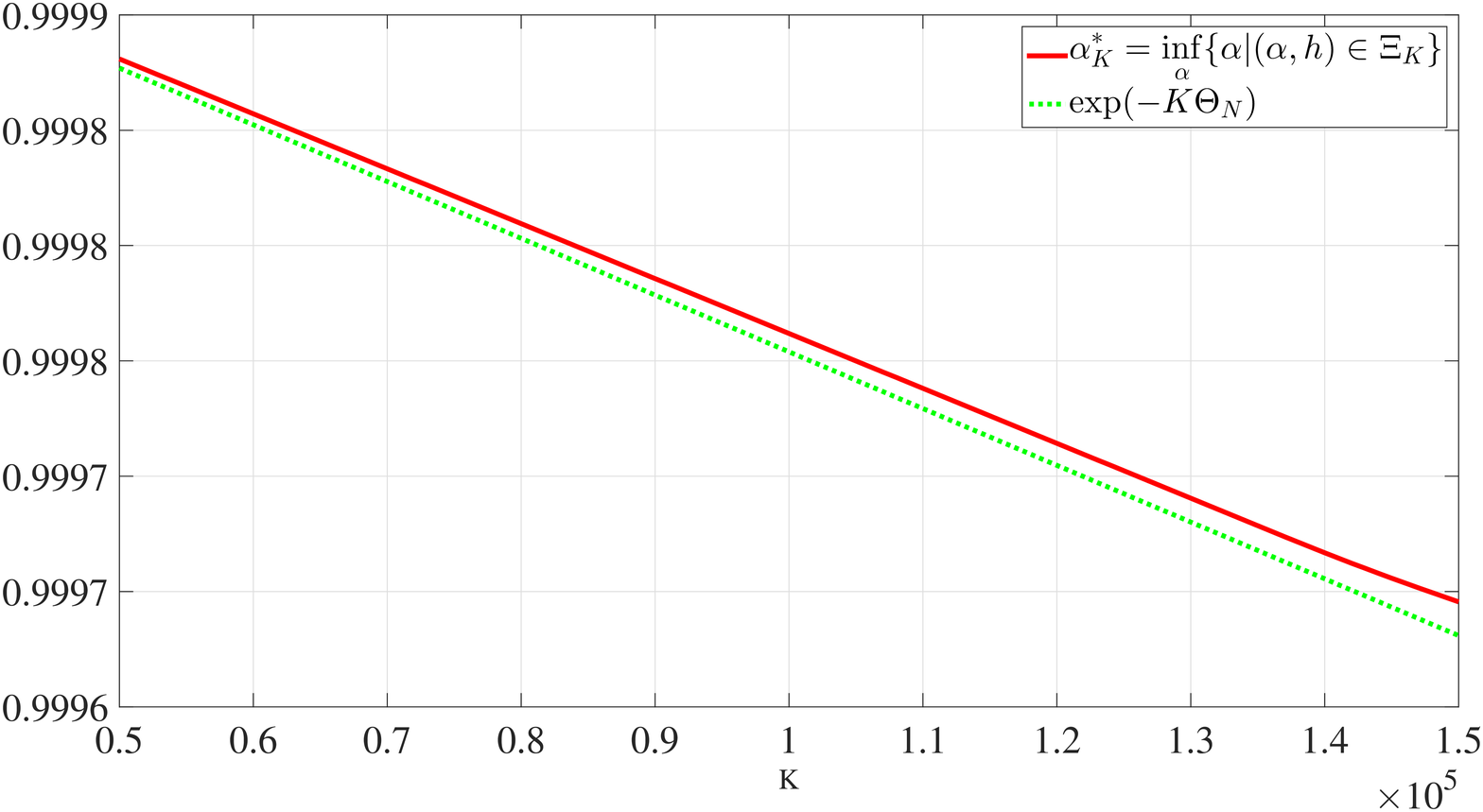}\\
  \caption{\small The minimal $\alpha_K^*$ and  $\exp(-K\Theta_N)$ with respect to  the data rate $K$ for $N=100$}\label{fig_exact_thm3}
\end{minipage}
\begin{minipage}[t]{0.04\textwidth}
\end{minipage}
\begin{minipage}[t]{0.46\textwidth}
\centering
    \includegraphics[width=3.6in]{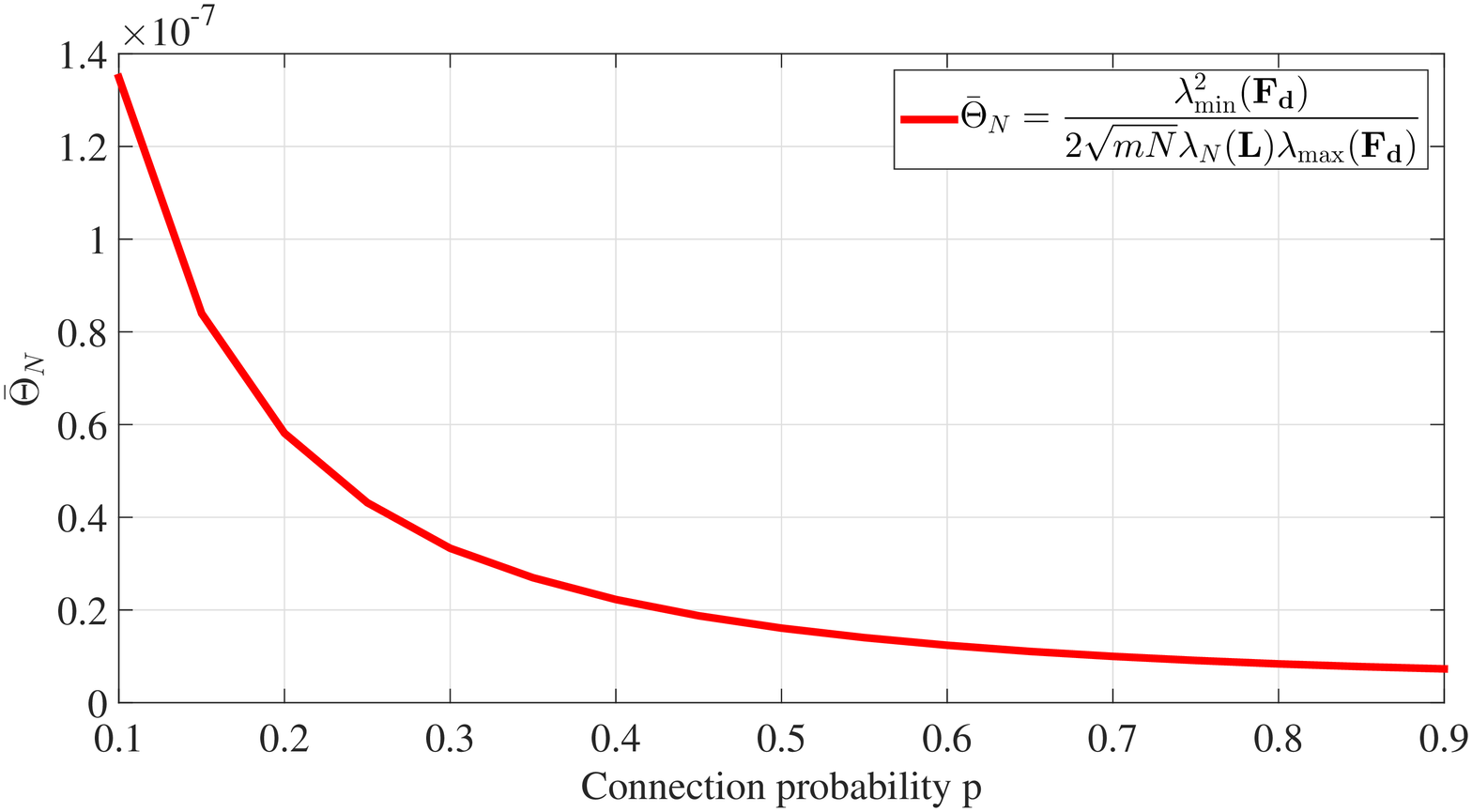}\\
  \caption{\small The mean $\bar{\Theta}_N$ for random graphs generated with different probability $p$.}\label{fig_exact_graphs}
\end{minipage}
\end{figure}

\noindent {\bf Example 2.} {\bf [Validation of Theorem \ref{thm-rate}]. } We let $N = 100$ and $m=5$. We randomly generate a matrix $\mathbf{H}$ and $\mathbf{z}$ such that $\mathbf{z}=\mathbf{H}\mathbf{y}$ has  a unique solution.
We set $\mathbf{L}$ as the Laplacian of a cycle  graph. Then  the constant  $\Theta_N= \frac{\lambda^2_{\min}(\mathbf{F_d}) }{ 2\sqrt{mN} \lambda_N(\mathbf{L}) \lambda_{\max}(\mathbf{F_d})}$ is fixed at $2.4910\times 10^{-9}$.
We let $K$  increase from $K=4\times 10^4$ to $K=1.5\times 10^5$ in steps of $1000$, and search for
the minimal $\alpha$ such that    $ (\alpha,h) \in \Xi_K $ for some $ h>0$ numerically for each $K$, i.e., $\alpha_K^*=\inf_{\alpha} \{\alpha| (\alpha,h)\in \Xi_K\}$.
Figure \ref{fig_exact_thm3} shows how $\alpha_K^*$ varies according to the data rate $K$,
 and implies that a higher data allows the selection of a smaller $\alpha$, and hence potentially leads to a faster convergence rate. Figure \ref{fig_exact_thm3} also  displays the trajectory of $\exp(-K\Theta_N)$ with respect to  $K$, and shows that
$\exp(-K\Theta_N)$ is quite close to $\alpha_K^*$ for  $N=100$, hence validates Theorem \ref{thm-rate}.
%

\noindent {\bf Example 3.}
Let $N=100$ and $m=10$. We randomly generate a matrix $\mathbf{H}$ and $\mathbf{z}$ such that $\mathbf{z}=\mathbf{H}\mathbf{y}$ has  a unique solution. It is easy to investigate how $\Theta_N$ depends on the network structure.
And for a complete graph, star graph and cycle graph, $\Theta_N$ takes values $6.9199\times 10^{-9}$, $1.8553\times 10^{-9}$, and $8.2899\times 10^{-8}$, respectively. This surprisingly indicates cycle graphs produce the fastest convergence compared to complete and star graphs.
We also compute $\Theta_N$ for Erd\H{o}s-R\`{e}nyi  random graphs $\mathcal{G}(N,p)$, where the possible connection between any two nodes is generated with a probability of $p$, independently of every other edge. We let $p$ increase from $0.1$ to $0.9$ in steps of $0.05$.
For each probability $p$, we  randomly generate  $10^3$ connected graphs with $\mathcal{G}(N,p)$, and compute the mean $\bar{\Theta}_N$. Figure \ref{fig_exact_graphs} shows how
$\bar{\Theta}_N$ varies along  with probability $p$, which decreases as the connection probability $p$ increases. This implies that $\alpha_K^*$, the fastest possible convergence rate under a fixed data rate $K$,   might increase with the increase of the connectivity of the graphs.
%

\subsection{ Discussion: Improve Robustness  with Damping}\label{subsec:exact:simulation:damping}

Convergence  of  Algorithm \ref{alg1}  relies on the equivalence
between  node $i$'s decoder output   $ \hat{\mathbf{x}}_{ij}(k) $ of its neighbor $j$' state and node $j$'s one-step prediction  $\mathbf{b}_j(k)$, which is characterized by   \eqref{equiv}. The theoretical and numerical results have shown the effectiveness of Algorithm \ref{alg1} when  \eqref{equiv} is satisfied.
In fact, \eqref{equiv} holds  when the encoder/decoder update \eqref{encoder}-\eqref{decoder} is exact and the following initialization condition is satisfied,
\begin{equation}\label{equ_assumption}
\mathbf{b}_j(0)=0,\; \hat{\mathbf{x}}_{ij}(0)=0, \; \forall j\in \mathcal{V},\forall i\in \mathcal{N}_j,
\end{equation}

However, there could exist initialization errors  in  \eqref{equ_assumption}. And due to round-off noises in the
 storage and manipulation of real-valued vectors  in digital computers, \eqref{encoder}-\eqref{decoder} may not be executed exactly.
With initialization errors in  \eqref{equ_assumption} and the round-off noises, the update of $\mathbf{b}_j(k)$, $\hat{\mathbf{x}}_{ij}(k)$ in encoder/decoder \eqref{encoder}-\eqref{decoder} is changed to
\begin{equation}
\begin{array}{ll}\label{equ_roundoff_encoder}
\mathbf{b}_j(0) &=\mathbf{I}^e_j,\;  \hat{\mathbf{x}}_{ij}(0)=\mathbf{I}^e_{ij}, \; \forall j\in \mathcal{V},\forall i\in \mathcal{N}_j,\\
\mathbf{b}_j(k)          & \triangleq   s(k-1)\mathbf{q}_j(k)  +  \mathbf{b}_j(k-1) + \varepsilon^b_j(k),\\
\hat{\mathbf{x}}_{ij}(k) & \triangleq   s(k-1)\mathbf{q}_j(k)  +  \hat{\mathbf{x}}_{ij}(k-1) + \varepsilon^x_{ij}(k).
\end{array}
\end{equation}
The  initialization errors  $\mathbf{I}^e_j, \mathbf{I}^e_{ij}$, and round-off noises $\varepsilon^b_j(k), \varepsilon^x_{ij}(k) $ will persist during the algorithm.

\begin{figure}[htbp]
\centering
\centering\includegraphics[width=4in]{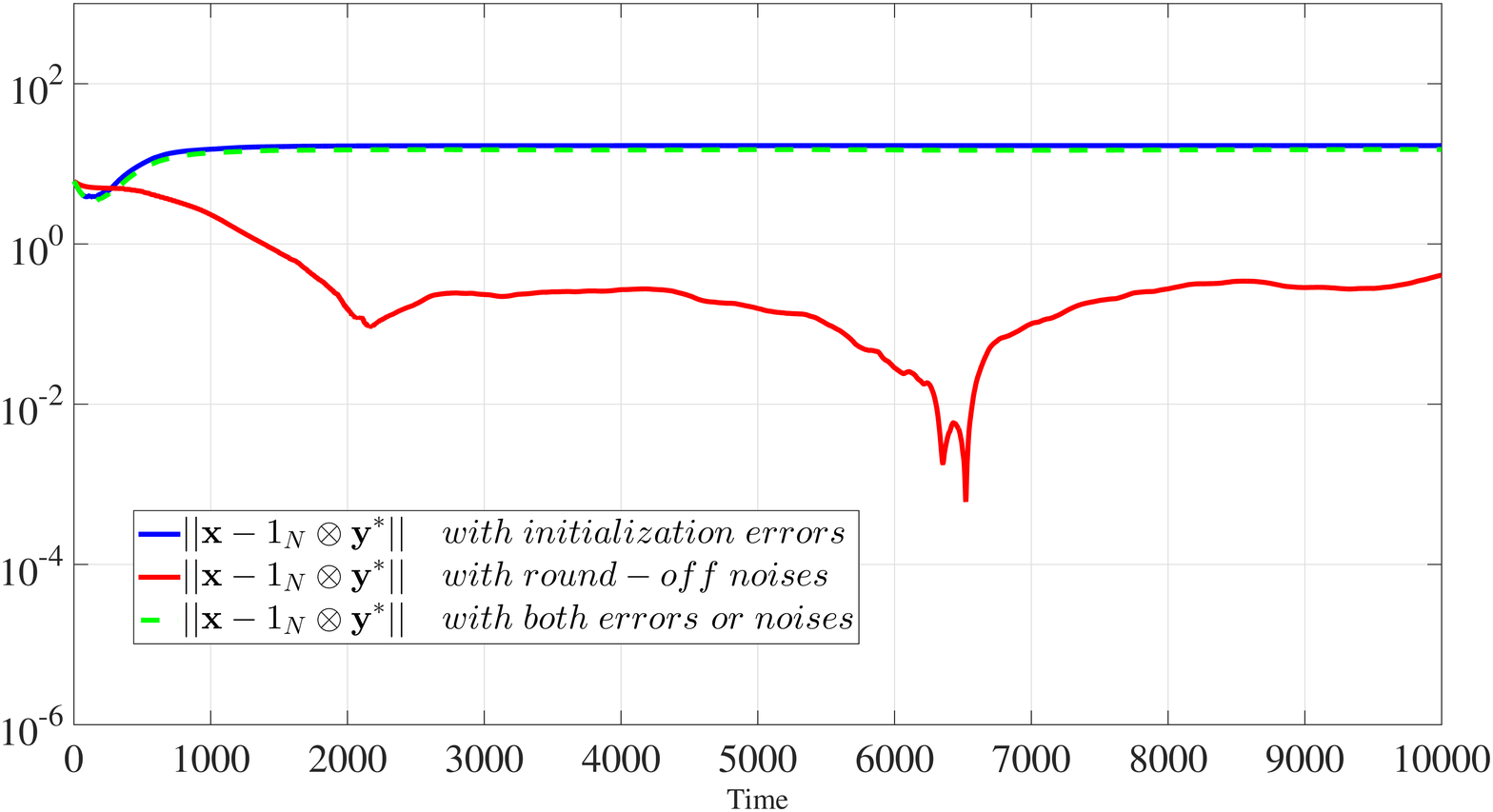}\\
  \caption{\small Trajectories of   $ ||\mathbf{x}-\mathbf{1}_N \otimes \mathbf{y}^* ||_2$ generated by  Algorithm \ref{alg1} for the following cases: (i) there exist initialization errors; (ii) there exist round-off noises; (iii) there exist both initialization errors and round-off noises.}\label{fig_noise_nodamping}
\end{figure}

\noindent{\bf Performance of Algorithm \ref{alg1} with  initialization errors and  round-off noises.} We continue to use the same $\mathbf{H}$ and $\mathbf{z}$ in \eqref{example_H_z}.
We set $h=0.0213$, $\alpha=0.998$, $K=300$ and $s(0)=10$. The initialization errors $\mathbf{I}^e_j$ and $\mathbf{I}^e_{ji}$ are independent and are randomly drawn from a uniform distribution on $[0,0.5]$, and
the round-off noises $\varepsilon^b_j(k)$, $\varepsilon^x_{ji}(k)$ are mutually independent random i.i.d. sequences with each value drawn from a uniform distribution on $[-1,1]\times 10^{-4}$. Figure \ref{fig_noise_nodamping} shows that Algorithm \ref{alg1} with \eqref{equ_roundoff_encoder} cannot ensure convergence when there exists initialization errors or round-off noises.
In fact, the error is very substantial in comparison to the average noise magnitude and the value of $||\mathbf{1}_N\otimes \mathbf{y}^*||_2$.


We propose to improve algorithm robustness by adding a damping term to encoder/decoder, where $\mathbf{b}_j(k)$ and $\hat{\mathbf{x}}_{ij}(k)$ are updated with
\begin{equation}
\begin{array}{ll}\label{equ_damping_encoder}
\mathbf{b}_j(0)&=\mathbf{I}^e_j,\;   \hat{\mathbf{x}}_{ij}(0)=\mathbf{I}^e_{ij}, \; \forall j\in \mathcal{V},\forall i\in \mathcal{N}_j,\\
\mathbf{b}_j(k)          & \triangleq   s(k-1)\mathbf{q}_j(k)  + \varrho  \mathbf{b}_j(k-1)+ \varepsilon^b_{j}(k),\\
\hat{\mathbf{x}}_{ij}(k) & \triangleq   s(k-1)\mathbf{q}_j(k)  + \varrho  \hat{\mathbf{x}}_{ij}(k-1)+ \varepsilon^x_{ij}(k),
\end{array}
\end{equation}
 where $\varrho\in (0,1)$ is a damping factor,  $\mathbf{I}^e_j, \mathbf{I}^e_{ij}$ are initialization errors, and $\varepsilon^b_j(k), \varepsilon^x_{ij}(k)$ are round-off noises.

\begin{figure}[htbp]
\centering
 \includegraphics[width=4in]{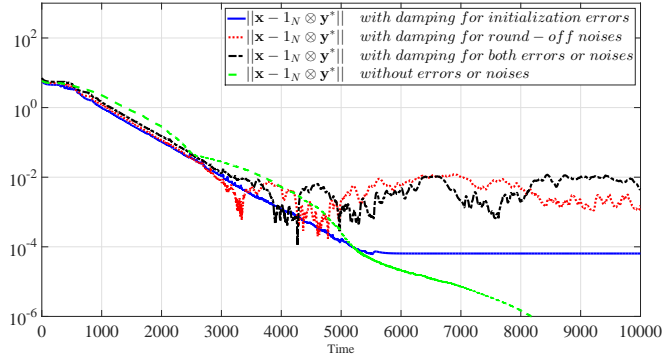}\\
  \caption{\small Trajectories of   $ ||\mathbf{x}-\mathbf{1}_N \otimes \mathbf{y}^* ||_2$ for (i) Algorithm \ref{alg1} with a damped encoder/decoder \eqref{equ_damping_encoder} with  initialization errors ; (ii) Algorithm \ref{alg1} with a damped encoder/decoder \eqref{equ_damping_encoder} with   round-off noises; (iii) Algorithm \ref{alg1} with a damped encoder/decoder \eqref{equ_damping_encoder} with   both initialization errors and round-off noises; (iv) Algorithm \ref{alg1} without  errors or   noises.}\label{fig_noise_damping}
\end{figure}

Now, we adopt the same setting as {\bf Example 1}.  We run Algorithm \ref{alg1} with \eqref{equ_damping_encoder} when there are initialization errors and round-off noises, and also run Algorithm \ref{alg1} with \eqref{encoder}-\eqref{decoder} where there are no initialization errors and round-off noises, both with the same algorithm parameters.
The damping factor is $\varrho=0.95$.
Figure \ref{fig_noise_damping} displays the simulation results, which shows  that
 (i) the damping can significantly reduce but not fully eliminate the affect of initialization errors in the final computed output (ii) the effect of round-off noises can be tolerated  in the sense that $\mathbf{x}_i(k)$ will converge to a neighborhood of the exact solution within   a distance of similar magnitude to  the round-off noises.

The formal   convergence analysis of  Algorithm \ref{alg1} with a damped encoder/decoder update \eqref{equ_damping_encoder} is  challenging because there will be a nonlinear coupling  between the damping factor $\varrho$ and all other   parameters, and  the errors and noises as well as $\varrho$ will enter the update equation of $\mathbf{x}(k)$ and $\mathbf{b}(k)$ in \eqref{compact-form-x}-\eqref{compact-form-xi}. Thereby, we leave the theoretical analysis of  \eqref{equ_damping_encoder} as a future research problem.

\subsection{Proof of Statements}\label{Sec1:proof}
\subsubsection{Preliminary Lemmas}
We first  give a reformulation of the  recursion for Algorithm \ref{alg1}.

\begin{lemma}\label{lem1}
Let {\bf A1} and {\bf A3} hold.
 Define
 \begin{align*}
 &   \mathbf{w}_i(k)=\mathbf{x}_i(k)-\mathbf{y}^*  ,  ~\mathbf{w}(k)=col\{  \mathbf{w}_1(k),\dots, \mathbf{w}_N(k)\}\\
 &   \mathbf{e}_i(k)=\mathbf{x}_i(k)- \mathbf{b}_i(k),  ~\mathbf{e}(k)=col\{  \mathbf{e}_1(k),\dots, \mathbf{e}_N(k)\},
  \bm{\omega}(k) \triangleq  {\mathbf{w}(k) \over s(k)},  ~ {\rm and}~  \bm{\varepsilon}(k) \triangleq {\mathbf{e}(k)\over s(k)}.
  \end{align*}   Then   the following hold:
 \begin{align}
    &  \bm{\omega}(k+1)=  \alpha^{-1}      \mathbf{P}_h   \bm{\omega}(k)
+   \alpha^{-1} h   \mathbf{L}\otimes \mathbf{I}_m \bm{\varepsilon}(k)  \label{cf-w} \\
    &  \bm{\varepsilon}(k+1) = \alpha^{-1} \big( \bm{\theta}(k)  - Q_K \left(\bm{\theta}(k)\right) \big) \label{cf-z},
\end{align}
where $  \mathbf{P}_h \triangleq  \mathbf{ I}_{mN}  - h  \mathbf{F_d}$ with $ \mathbf{F_d}= \mathbf{L}\otimes \mathbf{I}_m +  \mathbf{H_d}$,   and $\bm{\theta}(k)$ is defined as
 \begin{equation}  \label{compact-form2-xi2}
\begin{split}
\bm{\theta}(k)   \triangleq  & \left( \mathbf{I}_{mN}+ h  \mathbf{L}\otimes \mathbf{I}_m \right) \bm{\varepsilon}(k)  -
h \mathbf{F_d}\bm{\omega}(k).
\end{split}
\end{equation}
\end{lemma}
{\em Proof.}
Since $ j\in \mathcal{N}_i \Leftrightarrow  i\in \mathcal{N}_j $,  by using \eqref{equiv} and $ \mathbf{e}_i(k)=\mathbf{x}_i(k)- \mathbf{b}_i(k)$,  we have the following:

    \begin{equation}\label{consensus}
\begin{split}
 &\sum_{j\in \mathcal{N}_i}  \big(\hat{\mathbf{x}}_{ij}(k)-\mathbf{ b}_i(k)\big) = \sum_{j\in \mathcal{N}_i} \big( \mathbf{ b}_j(k)-\mathbf{ b}_i(k)\big) \\
 &= \sum_{j\in \mathcal{N}_i} \Big[\left(\mathbf{ x}_j(k)-\mathbf{ x}_i(k)\right)
 - \left(\mathbf{ x}_j(k)- \mathbf{ b}_j(k)\right)  
+\left(\mathbf{ x}_i(k)-\mathbf{b}_i(k)\right)\Big] \\
 &= \sum_{j\in \mathcal{N}_i}  \big[\left(\mathbf{ x}_j(k)-\mathbf{ x}_i(k)\right) - \left(\mathbf{ e}_j(k)-\mathbf{e}_i(k)\right)\big] .
\end{split}
\end{equation}
Recall that  $\mathbf{y}^*$ is the unique solution  to \eqref{LinearEquation} such that  $h_i^T \mathbf{y}^*=z_i~\forall i \in \mathcal{V} $.   Then by $  \mathbf{w}_i(k)=\mathbf{x}_i(k)-\mathbf{y}^* $, there holds
\begin{align*}
  \mathbf{h}_i\mathbf{h}_i^\top \mathbf{w}_i(k)
  &   = \mathbf{h}_i\mathbf{h}_i^\top \left( \mathbf{x}_i(k)-\mathbf{y}^*\right) = \mathbf{h}_i\mathbf{h}_i^\top   \mathbf{x}_i(k)-\mathbf{h}_iz_i  .
  \end{align*}
Also, using \eqref{encoder}, \eqref{quantized-do}, \eqref{consensus}, $\gamma(k)\equiv 1$, and the definition of $  \mathbf{H_d}$ in   \eqref{def}, leads to
  \begin{align}
\mathbf{x}(k+1) =&\mathbf{x}(k)- h  \mathbf{L}\otimes \mathbf{I}_m  \mathbf{x}(k) 
+ h  \mathbf{L}\otimes \mathbf{I}_m \mathbf{e}(k) -h   \mathbf{H_d}\mathbf{w}(k),\label{compact-form-x}
\\   \mathbf{b}(k+1) =&s(k)  Q_K\Big({\mathbf{x}(k+1)-\mathbf{b}(k)\over s(k)} \Big) +\mathbf{b}(k).\label{compact-form-xi}
\end{align}
Because $\mathbf{L} \mathbf{1}_N=\mathbf{0}_N,$ the following holds:
\begin{equation}\label{error-equiv1}
\begin{split}
&   \mathbf{L}\otimes \mathbf{I}_m  \mathbf{x}(k)  = \mathbf{L}\otimes \mathbf{I}_m  \mathbf{x}(k)-\left(\mathbf{L} \mathbf{1}_N\otimes \mathbf{I}_m\right)  \mathbf{y}^* \\&= \mathbf{L}\otimes \mathbf{I}_m \left( \mathbf{x}(k)-\mathbf{1}_N \otimes \mathbf{I}_m \mathbf{y}^*\right)= \mathbf{L}\otimes \mathbf{I}_m \mathbf{w}(k).
\end{split}
\end{equation}
 Now, by  subtracting $\mathbf{1}_N \otimes \mathbf{I}_m \mathbf{y}^*$ from both sides of   \eqref{compact-form-x}  and by substituting   \eqref{error-equiv1}, using   $\mathbf{w}(k)=\mathbf{x}(k )-\mathbf{1}_N \otimes \mathbf{I}_m \mathbf{y}^*$ and  $  \mathbf{F_d}= \mathbf{L}\otimes \mathbf{I}_m +  \mathbf{H_d}$, we obtain that
 \begin{align}
\mathbf{w}(k+1)&=\big(\mathbf{I}_{mN}- h \mathbf{F_d} \big)\mathbf{w}(k)
  + h  \mathbf{L}\otimes \mathbf{I}_m \mathbf{e}(k).\label{compact-form2-x}
\end{align}
Dividing  both sides of the above equation by $s(k+1)$, using  $s(k+1)=\alpha s(k) $ and   definitions of $      \mathbf{P}_h$,
 $\bm{\omega}(k) $ and  $\bm{\varepsilon}(k) $,  we obtain \eqref{cf-w}.

By subtracting  $\mathbf{ b}(k)$ from both sides of   \eqref{compact-form-x},   using \eqref{error-equiv1} and $\mathbf{e}(k)= \mathbf{x}(k )-\mathbf{ b}(k)$,  we obtain that
 \begin{align*}
\mathbf{x}(k+1)-\mathbf{ b}(k)&=\left( \mathbf{I}_{mN}+ h  \mathbf{L}\otimes \mathbf{I}_m \right) \mathbf{e}(k)
 -h \left(  \mathbf{L}\otimes \mathbf{I}_m +    \mathbf{H_d}\right) \mathbf{w}(k)  .
\end{align*}
Then   by using $\mathbf{w}(k)=s(k)\bm{\omega}(k)  $, $ \mathbf{e}(k)=s(k) \bm{\varepsilon}(k) ,$  $ \mathbf{F_d}= \mathbf{L}\otimes \mathbf{I}_m +  \mathbf{H_d}$,  and   \eqref{compact-form2-xi2},  we obtain that
 \begin{equation*}
\begin{split}
 \mathbf{x}(k+1)-\mathbf{ b}(k)   =    s(k) \bm{\theta}(k).
\end{split}
\end{equation*}
Now  recalling  that $\mathbf{ e}(k+1)=\mathbf{x}(k+1)- \mathbf{ b}(k+1)$  together with  \eqref{compact-form-xi},    the following holds:
\begin{equation*} \label{compact-form2-xi}
\begin{split}
\mathbf{ e}(k+1)&= \mathbf{x}(k+1)-\mathbf{ b}(k)-s(k)  Q_K\Big({\mathbf{x}(k+1)-\mathbf{ b}(k) \over s(k)}\Big) 
=s(k)\left(\bm{\theta}(k)- Q_K(\bm{\theta}(k))\right).
\end{split}
\end{equation*}
Dividing  both sides of the above  equation  by  $s(k+1)$ and using   $s(k+1)=\alpha s(k) ,$ we obtain \eqref{cf-z}.
\hfill$\square$

  \subsubsection{Proof of Proposition  \ref{prp1} }\label{proof-thm01}
 The proof of  non-saturation of the uniform   quantizer is equivalent to showing that for any $k\geq 0$, $ \bm{\theta}(k)$ defined  by \eqref{compact-form2-xi2} satisfies  $\| \bm{\theta}(k)\|_{\infty}<  K+{1\over 2}$.
 The proof  of Proposition  \ref{prp1}  will use induction, and we begin by showing the quantizer  is not saturated at $k=0.$

By using $\mathbf{ b}_i(0)=0 ~ \forall i\in V,$ we obtain that $\mathbf{ e}(0)= \mathbf{ x}(0)$ and  $\bm{\varepsilon}(0)= {\mathbf{ x}(0)/ s(0)}.$ Then by {\bf A2}, we have
\begin{align}\label{bd-z0}
 \| \bm{\varepsilon}(0) \|_{\infty}={ \|  \mathbf{ x}(0) \|_{\infty} \over s(0) } \leq {C_x / s(0)}.
\end{align}
By \eqref{error-equiv1},   and by recalling  that $\bm{\omega}(0)={ \mathbf{w}(0)/ s(0)}$ and  $\bm{\varepsilon}(0)={\mathbf{x}(0) / s(0)}$, we obtain that
    $$   \mathbf{L}\otimes \mathbf{I}_m \bm{\omega}(0)=   {   \mathbf{L}\otimes \mathbf{I}_m \mathbf{x}(0) / s(0)}= \mathbf{L}\otimes \mathbf{I}_m \bm{\varepsilon}(0).$$
 Then  by definition   \eqref{compact-form2-xi2}
there holds
\begin{align*}
\bm{\theta}(0) & =  \left( \mathbf{I}_{mN}+ h  \mathbf{L}\otimes \mathbf{I}_m \right) \bm{\varepsilon}(0)-h\left(  \mathbf{L}\otimes \mathbf{I}_m +  \mathbf{H_d}\right) \bm{\omega}(0)
=\bm{\varepsilon}(0) - h  \mathbf{H_d}{ \mathbf{w}(0) / s(0)}.
\end{align*}
As a result, by {\bf A2},  \eqref{def-s0}  and \eqref{bd-z0}  we have the following:
\begin{align*}
 \| \bm{\theta}(0)\|_{\infty}   &\leq \| \bm{\varepsilon}(0) \|_{\infty} +  {h  \|   \mathbf{H_d}  \mathbf{w}(0) \|_{\infty} /s(0)}
\leq  \left(C_x+h  \|   \mathbf{H_d}  \|_{\infty} C_w \right)/s(0)<  K+{1\over 2}.
\end{align*}
Hence, when $k=0,$ the quantizer is unsaturated.
Now for the induction, we assume  that when $k=0,\dots, p$, the quantizer is not saturated.
Then by  \eqref{cf-z} we have that
\begin{align}\label{math-induct-z}
\sup_{1\leq k \leq p+1}  \| \bm{\varepsilon}(k)  \|_{\infty} \leq {1 \over 2 \alpha}.
\end{align}
We proceed  to show that the quantizer is unsaturated for  $k=p+1$.

From \eqref{cf-w} it follows that
\begin{equation}\label{recursion-w}
\begin{split}
    \bm{\omega}(p+1)  &= \left(  \alpha^{-1}     \mathbf{P}_h  \right)^{p+1} \bm{\omega}(0)
    +\alpha^{-1} h  \left(  \alpha^{-1}     \mathbf{P}_h  \right)^{p }   \mathbf{L}\otimes \mathbf{I}_m \bm{\varepsilon}(0)    +\alpha^{-1} h \sum_{i=0}^{p-1} \left(  \alpha^{-1}     \mathbf{P}_h  \right)^{i }   \mathbf{L}\otimes \mathbf{I}_m \bm{\varepsilon}(p-i).
\end{split}
\end{equation}
We now estimate the three terms on the right-hand side of the above equation separately.
Note that  any given $h>0$,    the eigenvalues of $ \mathbf{P}_h=\mathbf{ I}_{mN}  - h  \mathbf{F_d}$  are sorted in an ascending order as  $1-h  \lambda_{\max}(\mathbf{F_d}) \leq \dots \leq1-h  \lambda_{\min}(\mathbf{F_d}),$ and  there exists  a unitary matrix $\mathbf{U}$ such that $ \mathbf{U}^T   \mathbf{P}_h \mathbf{U}= {\rm diag}  \big\{1-h  \lambda_{\max}(\mathbf{F_d}), \dots ,1-h  \lambda_{\min}(\mathbf{F_d}) \big\}\triangleq \bm{\Lambda}.$
Therefore,
\begin{align}\label{pk}
\left(   \mathbf{P}_h \right)^k =\left(\mathbf{U}  \bm{\Lambda} \mathbf{U}^T\right)^k=  \mathbf{U} \bm{\Lambda}^k \mathbf{U}^T.
\end{align}
By using the definition of  $\bm{\Lambda}$,  we obtain that
$$\|\bm{\Lambda}\|_2 =\max \left \{ |1-h  \lambda_{\min}(\mathbf{F_d}) |,  |1-h  \lambda_{\max}(\mathbf{F_d})|\right\}.$$
Thus, by using   $h \in \big(0, {2\over \lambda_{\min}(\mathbf{F_d})+\lambda_{\max}(\mathbf{F_d})}\big)$ and
    \cite[Lemma 3.1]{litao2011},  there holds $\|\bm{\Lambda}\|_2=1-h \lambda_{\min}(\mathbf{F_d})= \rho_h $.
 For the first term,  using \eqref{pk}, $\|\mathbf{U}\|_2=1$ and
$\|\mathbf{ x}\|_{\infty} \leq \|\mathbf{ x}\|_2\leq \sqrt{m} \|\mathbf{ x}\|_{\infty}$ for any $\mathbf{ x}\in \mathbb{R}^m, $ we have
\begin{align}\label{recursion-term1}
  & \| \left(  \alpha^{-1}     \mathbf{P}_h  \right)^{p+1} \bm{\omega}(0) \|_2
  \leq   \left \| \mathbf{U} \left(  {  \bm{\Lambda} / \alpha}  \right)^{p+1} \mathbf{U}^T  \right \|_2  \|\bm{\omega}(0) \|_2 \notag
\\&\leq  \left(  { \rho_h \over \alpha}  \right)^{p+1} {  \|\mathbf{w}(0) \|_2 \over s(0)}
 \leq  \left(  { \rho_h \over \alpha}  \right)^{p+1} { \sqrt{mN} \|\mathbf{w}(0) \|_{\infty} \over s(0)}
 <    { \sqrt{mN}C_w  \over s(0)} \left(  { \rho_h \over \alpha}  \right)^{p +1}
~{\scriptstyle \left(  \mathrm{by ~{\bf A2}}  \right) }.
\end{align}
For the second term of \eqref{recursion-w},  using \eqref{pk},  \eqref{bd-z0}, and $ \| \mathbf{L} \|_2 =\lambda_N(\mathbf{L})$   we obtain the following:
\begin{equation}\label{recursion-term2}
\begin{split}
&  \|  \alpha^{-1} h  \left(  \alpha^{-1}     \mathbf{P}_h  \right)^{p }   \mathbf{L}\otimes \mathbf{I}_m \bm{\varepsilon}(0)   \|_2 \\& \leq
 \alpha^{-1} h   \left \| \mathbf{U} \left(  {  \bm{\Lambda}/ \alpha}  \right)^{p } \mathbf{U}^T  \right \|_2  \| \mathbf{L} \|_2 \| \bm{\varepsilon}(0)   \|_2
  \leq   {hC_x\sqrt{mN}  \over   \alpha s(0)}   \lambda_N(\mathbf{L})  \left(  { \rho_h \over \alpha}  \right)^{p }
\end{split}
\end{equation}
Similarly, for the last term of \eqref{recursion-w}, by
\begin{align*}
& \Big \|  \sum_{i=0}^{p-1} \left(  \alpha^{-1}     \mathbf{P}_h  \right)^{i }  \Big  \|_2 \leq\sum_{i=0}^{p-1} \left  \|   \left(  \alpha^{-1}     \mathbf{P}_h  \right)^{i }  \right  \|_2
  \leq\sum_{i=0}^{p-1} \left(  { \rho_h \over \alpha}  \right)^i={1-\left(  { \rho_h/ \alpha}  \right)^p \over  1-  { \rho_h/ \alpha}  },
\end{align*}
and  by \eqref{math-induct-z}  we have that
\begin{align}\label{recursion-term3}
& \Big\|   \alpha^{-1} h \sum_{i=0}^{p-1} \left(  \alpha^{-1}     \mathbf{P}_h  \right)^{i }   \mathbf{L}\otimes \mathbf{I}_m \bm{\varepsilon}(p-i)\Big\|_{2}
\notag  \\ &  \leq      {h \sqrt{mN} \over  \alpha }     \| \mathbf{L} \|_2\sup_{1\leq k \leq p+1}  \| \bm{\varepsilon}(k)  \|_{\infty} \Big  \|  \sum_{i=0}^{p-1} \left(  \alpha^{-1}     \mathbf{P}_h  \right)^{i }  \Big  \|_2
 \leq      {h \sqrt{mN}  \lambda_N(\mathbf{L})  \over 2 \alpha(\alpha-\rho_h)}
\Big (1-\left(  { \rho_h \over \alpha}  \right)^p  \Big  ).
\end{align}

Since $\alpha \in ( \rho_h,1)$, by using \eqref{recursion-w} and \eqref{recursion-term1}-\eqref{recursion-term3},  we have that
\begin{equation}\label{recursion-w2}
\begin{split}
 & \|  \bm{\omega}(p+1) \|_{\infty} \leq \|  \bm{\omega}(p+1) \|_2 \\& \leq
{ \sqrt{mN} \over  \alpha} \max\left \{  { \rho_h  C_w  + hC_x  \lambda_N(\mathbf{L})   \over   s(0)}  ,   {h  \lambda_N(\mathbf{L}) \over 2(\alpha-\rho_h)} \right \}  \leq {h \sqrt{mN} \lambda_N(\mathbf{L}) \over 2 \alpha(\alpha-\rho_h)},
\end{split}
\end{equation}
where the last inequality follows by \eqref{def-s0}.
This together with $\|\mathbf{ x}\|_{\infty} \leq \|\mathbf{ x}\|_2$,
\eqref{def-M}, \eqref{compact-form2-xi2},   and \eqref{math-induct-z}  leads to
\begin{equation*}
\begin{split}
&  \| \bm{\theta}(p+1)\|_{\infty}
 \leq  \|  \left( \mathbf{I}_{mN} + h  \mathbf{L}\otimes \mathbf{I}_m \right) \bm{\varepsilon}(p+1) \|_{\infty}+ \| h   \mathbf{F_d} \bm{\omega}(p+1)\|_{\infty}
\\& \leq \|    \mathbf{I}_{mN}+ h  \mathbf{L}\otimes \mathbf{I}_m  \|_{\infty}\| \bm{\varepsilon}(p+1) \|_{\infty}
+ h\|    \mathbf{F_d} \|_2 \| \bm{\omega}(p+1)\|_2 \\& \leq {1+2hd^* \over 2 \alpha}+
{h^2 \sqrt{mN}  \lambda_N(\mathbf{L}) \lambda_{\max} (\mathbf{F_d})\over 2 \alpha(\alpha-\rho_h)}
= M(\alpha, h)  \leq\mathcal{K}(\alpha, h )+{1 \over 2} \leq K+{1\over 2} .
\end{split}
\end{equation*}
As a result, when $k=p+1$, the quantizer is also unsaturated.
Therefore, by induction, we conclude that if a $(2K+1)$-levels uniform quantizer is applied, then the  quantizer will   never be saturated.
\hfill $\blacksquare$

\subsubsection{Proof of Theorem \ref{thm1}}\label{proof-thm1}
 Since the conditions  required by Proposition  \ref{prp1} are the same as those used in Theorem  \ref{thm1}, the  quantizer will   never be saturated by   Proposition  \ref{prp1}.
Then by     \eqref{cf-z}  we conclude  that  $\sup_{k\geq 1}  \| \bm{\varepsilon}(k)  \|_{\infty} \leq {1 / 2 \alpha},$
and hence  \eqref{recursion-w2}  holds for any $p\geq 0.$
Thus,
\begin{equation*}
 \limsup_{k\to \infty}  \|  \bm{\omega}(k) \|_2 \leq  {h \sqrt{mN}  \lambda_N(\mathbf{L}) \over 2 \alpha(\alpha-\rho_h)}.
\end{equation*}
 Then by  using $\mathbf{w}(k)=s(0) \alpha^k \bm{\omega}(k) $ and  $\mathbf{w}(k)=\mathbf{x}(k)- \mathbf{1}_N \otimes \mathbf{y}^*$,  we obtain \eqref{rate} and  \eqref{convergence}.
\hfill $\blacksquare$

\subsubsection{Proof of Theorem \ref{thm2}}\label{proof-thm2}
 \noindent (i) By using   $\rho_h =1-h\lambda_{\min}(\mathbf{F_d})$  and  $M(\alpha, h )$ defined in \eqref{def-M}, there holds:
  \begin{align}\label{def-M2}
M(\alpha, h )={ 1+2hd^* \over 2 \alpha}+
{h^2 \sqrt{mN} \lambda_N(\mathbf{L}) \lambda_{\max}(\mathbf{F_d})\over 2 \alpha\big(\alpha-(1-h\lambda_{\min}(\mathbf{F_d}))\big)} .
\end{align}
Noting that
$$\lim_{h\to 0}{ 1+2hd^* \over 2 }+
{h  \sqrt{mN} \lambda_N(\mathbf{L}) \lambda_{\max}(\mathbf{F_d})\over 2   \lambda_{\min}(\mathbf{F_d}) }={1\over 2},   $$
then for any given $K\geq 1$ there exists $h^* \in \left(0, {2\over \lambda_{\min}(\mathbf{F_d})+\lambda_{\max}(\mathbf{F_d})}\right)$ such that
$$ { 1+2h^*d^* \over 2 }+
{h^*  \sqrt{mN} \lambda_N(\mathbf{L}) \lambda_{\max}(\mathbf{F_d})\over 2   \lambda_{\min}(\mathbf{F_d}) }\leq K .  $$
By \eqref{def-M2} it follows that
\begin{align*}
& \lim_{\alpha \to 1}M(\alpha, h^*) ={ 1+2h^*d^* \over 2 } +
{h^*  \sqrt{mN} \lambda_N(\mathbf{L}) \lambda_{\max}(\mathbf{F_d})\over 2   \lambda_{\min}(\mathbf{F_d}) }\leq K,\end{align*}
and hence there exists $\alpha^* \in (1-h^*\lambda_{\min}(\mathbf{F_d}),1)$ such that $M(\alpha^*, h^*)  <  K+{1\over 2}.$ Thus,   $(\alpha^*, h^*) \in  \Xi_K  , $  and  hence $ \Xi_K  $ is nonempty.

\noindent (ii) For any $(\alpha,h) \in  \Xi_K $,  from  \eqref{def_omega} it follows that
$  h \in \left(0, {2\over \lambda_{\min}(\mathbf{F_d})+\lambda_{\max}(\mathbf{F_d})}\right) ,\alpha\in (1-h \lambda_{\min}(\mathbf{F_d}) ,1), $ and $M(\alpha, h)< K +{1\over 2}$. Then   $ \mathcal{K}(\alpha, h )\triangleq \left \lceil M(\alpha, h -{1\over 2}\right \rceil  \leq K $ together with  Theorem \ref{thm1} leads to the result (ii).
\hfill $\blacksquare$

\subsubsection{  Proof of Proposition \ref{lem-rate} }\label{proof-lem}
We first prove  $\bigcup_{\epsilon \in (0,1)} \Xi_{K,\epsilon}\subset \Xi_K$.
For any given $K\geq 1 $ and $ \epsilon \in (0,1)$, let $(\alpha,h)\in  \Xi_{K,\epsilon}$.
Then $\alpha-\rho_h= \epsilon h \lambda_{\min}(\mathbf{F_d})  >0$ by $\rho_h=1-  h \lambda_{\min}(\mathbf{F_d}) $, and
$\alpha \in (1-  h \lambda_{\min}(\mathbf{F_d}) ,1).$  Also,   by  the definition $ M(\alpha, h)$ in \eqref{def-M}, we obtain the following:
\begin{equation}\label{def-s-epsilon}
\begin{split}
 &S(\epsilon,h) \triangleq M(\alpha, h)={ 1+2hd^* \over 2 \left(1-(1-\epsilon) h \lambda_{\min}(\mathbf{F_d}) \right)}
 + {h  \sqrt{mN} \lambda_N(\mathbf{L}) \lambda_{\max}(\mathbf{F_d})\over 2 \epsilon   \lambda_{\min}(\mathbf{F_d})    \left(1-(1-\epsilon) h \lambda_{\min}(\mathbf{F_d}) \right)}
\\& ={\epsilon  \lambda_{\min}(\mathbf{F_d}) ( 1+2hd^*)    +h \sqrt{mN} \lambda_N(\mathbf{L}) \lambda_{\max}(\mathbf{F_d})\over 2\epsilon  \lambda_{\min}(\mathbf{F_d})  \left(1-(1-\epsilon) h \lambda_{\min}(\mathbf{F_d}) \right)  } .
\end{split}
\end{equation}
Then by using the definition of $ \hat{h}_{K,\epsilon} $  in \eqref{def-hath} and $h<\hat{h}_{K,\epsilon}$,  there holds  $ M(\alpha, h ) < K+{1\over 2}$. It is clear that for any $(\alpha,h)\in  \Xi_{K,\epsilon}$,  $h \in \left(0, {2\over \lambda_{\min}(\mathbf{F_d})+\lambda_{\max}(\mathbf{F_d})}\right)$. In summary, we have shown that  for any given $K\geq 1 $ and $ \epsilon \in (0,1)$,   $\Xi_{K,\epsilon}\subset \Xi_K.$ Thus, $\bigcup_{\epsilon \in (0,1)} \Xi_{K,\epsilon}\subset \Xi_K$.

We now validate  $\Xi_K\subset \bigcup_{\epsilon \in (0,1)} \Xi_{K,\epsilon}$.
For any  $(\alpha_0,h_0)\in \Xi_K$, by \eqref{def_omega} we have $h_0 \in \left(0, {2\over \lambda_{\min}(\mathbf{F_d})+\lambda_{\max}(\mathbf{F_d})}\right)$ and $\rho_{h_0} =1-h_0 \lambda_{\min}(\mathbf{F_d}) $.
Note that  $\alpha_0=1-(1-\epsilon_0) h_0 \lambda_{\min}(\mathbf{F_d}) \in (\rho_{h_0},1)$ with $\epsilon_0=1-{1-\alpha_0 \over  h_0 \lambda_{\min}(\mathbf{F_d}) }$. Then $\epsilon_0\in (0,1)$ and
$M(\alpha_0, h_0) = S(\epsilon_0,h_0) $, where  $S(\epsilon_0,h_0) $ is given by \eqref{def-s-epsilon} with $(\epsilon,h)$ replaced by $(\epsilon_0,h_0) $. This together with $M(\alpha_0, h_0) <K+{1\over 2}$  leads to $S(\epsilon_0,h_0) <K+{1\over 2}.$ This is equivalent to
\begin{align*}
& h_0<\hat{h}_{K,\epsilon_0 }= 2K  \epsilon_0  \lambda_{\min}(\mathbf{F_d}) \Big(\sqrt{mN} \lambda_N(\mathbf{L}) \lambda_{\max}(\mathbf{F_d})\\&+2 \epsilon_0  \lambda_{\min}(\mathbf{F_d})  d^*+ \epsilon_0  (1-\epsilon_0) (2K+1)    \lambda^2_{\min}(\mathbf{F_d})   \Big)^{-1}.
\end{align*}
Then by  $h_0 \in \left(0, {2\over \lambda_{\min}(\mathbf{F_d})+\lambda_{\max}(\mathbf{F_d})}\right)$,  we conclude that
$h_0\in (0,h_{K,\epsilon_0}^*)$. This together with $\epsilon_0\in (0,1)$  implies that  $(\alpha_0,h_0) \in  \Xi_{K,\epsilon_0}, $ and hence $\Xi_K\subset \bigcup_{\epsilon \in (0,1)} \Xi_{K,\epsilon}$.
This  completes the proof of Lemma \ref{lem-rate}.
\hfill $\blacksquare$

\subsubsection {Proof of  Theorem \ref{thm-rate}} \label{proof-thm-rate} For any given $K\geq 1 $,  define
  \begin{equation} \label{def-gammak}
\begin{split}  \Gamma_{K} \triangleq    \{ \alpha:  &\alpha=1-(1-\epsilon) h \lambda_{\min}(\mathbf{F_d}) , \epsilon \in (0,1), h\in (0,h_{K,\epsilon}^*)   \} .
\end{split}
\end{equation}
By $h_{K,\epsilon}^* \leq \hat{h}_{K,\epsilon} $ and \eqref{def-hath},  we know for any $  h\in (0,h_{K,\epsilon}^*) ,$
 $ h \leq 2K  \epsilon  \lambda_{\min}(\mathbf{F_d}) \left (  \sqrt{mN} \lambda_N(\mathbf{L}) \lambda_{\max}(\mathbf{F_d})  \right )^{-1}.$
Then  for any $\alpha \in \Gamma_K$ with $\epsilon \in (0,1), h\in (0,h_{K,\epsilon}^*),$  from $ (1-\epsilon)      \epsilon \leq {1\over 4}~ \forall \epsilon \in (0,1)$ it follows that
\begin{align*}
\alpha&=1-(1-\epsilon) h \lambda_{\min}(\mathbf{F_d}) >1- { 2K  (1-\epsilon)      \epsilon  \lambda^2_{\min}(\mathbf{F_d}) \over   \sqrt{mN} \lambda_N(\mathbf{L}) \lambda_{\max}(\mathbf{F_d}) }  \geq 1- { K  \lambda^2_{\min}(\mathbf{F_d}) \over 2\sqrt{mN} \lambda_N(\mathbf{L}) \lambda_{\max}(\mathbf{F_d})  },
\end{align*}
Thus,  the following holds for fixed $K$:
\begin{align*}
& { \inf_{\alpha  \in \Gamma_K}\alpha  \over  {\rm exp} \left( -{ K  \lambda^2_{\min}(\mathbf{F_d}) \over 2\sqrt{mN} \lambda_N(\mathbf{L}) \lambda_{\max}(\mathbf{F_d})  }\right)}   \geq  { 1- { K  \lambda^2_{\min}(\mathbf{F_d}) \over 2\sqrt{mN} \lambda_N(\mathbf{L}) \lambda_{\max}(\mathbf{F_d})  } \over {\rm exp} \left( -{ K  \lambda^2_{\min}(\mathbf{F_d}) \over 2\sqrt{mN} \lambda_N(\mathbf{L}) \lambda_{\max}(\mathbf{F_d})  }\right)},
\end{align*}
which together with  $\lim_{x\downarrow 0} {1-x \over  {\rm exp} (-x) }=0$ produces
\begin{equation}\label{bd1-alpha}
\begin{split}
\liminf_{N\to \infty }{ \inf_{\alpha  \in \Gamma_K}\alpha  \over  {\rm exp} \left( -{ K  \lambda^2_{\min}(\mathbf{F_d}) \over 2\sqrt{mN} \lambda_N(\mathbf{L}) \lambda_{\max}(\mathbf{F_d})  }\right)}  \geq  1.
\end{split}
\end{equation}

From \eqref{def-hath} it follows that  \begin{align*} \hat{h}_{K,\epsilon} \geq   &2K  \epsilon  \lambda_{\min}(\mathbf{F_d}) \Big(  \sqrt{mN} \lambda_N(\mathbf{L}) \lambda_{\max}(\mathbf{F_d})
+ 2    \lambda_{\min}(\mathbf{F_d})  d^*+  (2K+1)    \lambda^2_{\min}(\mathbf{F_d})   \Big )^{-1}.
\end{align*}
This together with $$\inf_{h\in (0,h_{K,\epsilon}^*)}~ \alpha \leq 1-(1-\epsilon) \hat{h}_{K,\epsilon}   \lambda_{\min}(\mathbf{F_d}) $$ implies
\begin{align*}
&\inf_{ \alpha \in \Gamma_K}\alpha  \leq    1-  { \max_{\epsilon \in (0,1)}~ 2\epsilon (1-\epsilon) K  \lambda^2_{\min}(\mathbf{F_d}) \over 2    \lambda_{\min}(\mathbf{F_d})  d^*+  \sqrt{mN} \lambda_N(\mathbf{L}) \lambda_{\max}(\mathbf{F_d})+  (2K+1)    \lambda^2_{\min}(\mathbf{F_d})  }\\& = 1-
{ K  \lambda^2_{\min}(\mathbf{F_d}) /2\over
2 \lambda_{\min}(\mathbf{F_d})  d^*   + \sqrt{mN} \lambda_N(\mathbf{L}) \lambda_{\max}(\mathbf{F_d}) +(2K+1)    \lambda^2_{\min}(\mathbf{F_d}) }
\end{align*}
Then we have that for any given $K\geq 1,$
\begin{align*}
&  { \inf_{\alpha  \in \Gamma_K}\alpha  \over  {\rm exp} \left( -{ K  \lambda^2_{\min}(\mathbf{F_d}) \over 2\sqrt{mN} \lambda_N(\mathbf{L}) \lambda_{\max}(\mathbf{F_d})  }\right)} \\ &\leq   { 1- { K  \lambda^2_{\min}(\mathbf{F_d}) \over 2\sqrt{mN} \lambda_N(\mathbf{L}) \lambda_{\max}(\mathbf{F_d})  } \over {\rm exp} \left( -{ K  \lambda^2_{\min}(\mathbf{F_d}) \over 2\sqrt{mN} \lambda_N(\mathbf{L}) \lambda_{\max}(\mathbf{F_d})  }\right)}
 \times {  1-{ K  \lambda^2_{\min}(\mathbf{F_d})/2 \over  \sqrt{mN} \lambda_N(\mathbf{L}) \lambda_{\max}(\mathbf{F_d}) +
2  \lambda_{\min}(\mathbf{F_d})  d^*   + (2K+1)    \lambda^2_{\min}(\mathbf{F_d}) } \over  1- { K  \lambda^2_{\min}(\mathbf{F_d}) /2\over  \sqrt{mN} \lambda_N(\mathbf{L}) \lambda_{\max}(\mathbf{F_d})  } },
\end{align*}
which together with
$\lim_{N \to \infty} {1-c_1/(\sqrt{N}+c_2) \over 1-c_1/\sqrt{N}}=1$ and  $\lim_{x\downarrow 0} {1-x \over  {\rm exp} (-x) }=0$    gives
\begin{align}\label{bd2-alpha}
\limsup_{N\to \infty }{ \inf_{\alpha  \in \Gamma_K}\alpha  \over  {\rm exp} \left( -{ K  \lambda^2_{\min}(\mathbf{F_d}) \over 2\sqrt{mN} \lambda_N(\mathbf{L}) \lambda_{\max}(\mathbf{F_d})  }\right)}  \leq  1.
\end{align}
Using Lemma \ref{lem-rate} and \eqref{def-gammak}, we have the following:
$$\inf_{ (\alpha,h)  \in \Xi_K}=  \inf_{ (\alpha,h)  \in \cup_{\epsilon \in (0,1)} \Xi_{K,\epsilon}}=\inf_{\epsilon \in (0,1)}  \inf_{ (\alpha,h)  \in \Xi_{K,\epsilon}}=\inf_{  \alpha    \in \Gamma_K}.$$
By this, using   \eqref{bd1-alpha}  and  \eqref{bd2-alpha}, we obtain  \eqref{bd-alpha}.
\hfill $\blacksquare$

\section{ Least-Squares Solver }\label{sec:leastsquares}
In this section, we investigate  the case  ${\rm rank}(\mathbf{H})=m$ and $\mathbf{z}\notin {\rm span}(\mathbf{H})$.
Then  equation \eqref{LinearEquation} does not have exact solutions,  while  a  least-squares
solution  is defined as the solution to the optimization problem \eqref{LS}.
We consider Algorithm \ref{alg1},  then show the convergence results regarding the  quantization level  along with
the data rate analysis, and demonstrate the results with numerical simulations.

  \subsection{Convergence Results}
 Assumptions {\bf A1}, {\bf A2} and {\bf A3} are no longer in force, instead, we impose   the following conditions on the initial states and step-size.

\noindent{\bf A4}  ${\rm rank}(\mathbf{H})=m$ and $\mathbf{z}\notin {\rm span}(\mathbf{H})$.

\noindent  {\bf A5}  $\max_i \|\mathbf{x}_i(0)\|_{\infty}\leq C_x$ for  constant  $C_x>0$.

\noindent  {\bf A6}  (i)  $\gamma(0)= 1$, $\gamma(k) \downarrow 0,$ $\sum_{k=1}^{\infty} \gamma(k)=\infty,$  (ii) $s(k) = s_r \gamma(k)  $ for some $s_r>0$,  and (iii)  $1<\beta(k+1)<\beta(k)$  for any $k\geq 0$, where $\beta(k) \triangleq  {\gamma(k) \over \gamma(k+1)}$.

\begin{remark} We  now specify how to choose $\gamma(k)$ to make {\bf A6} hold. Set $\gamma(k)={k_0^{\delta}\over (k+k_0)^{\delta}}$ for some $\delta\in({1\over 2},1] $, where $k_0={1 \over \beta(0)^{1/\delta}-1}$.  Then it  is seen that  $\gamma(0)= 1, \gamma(k)  \downarrow 0 $ and $\sum_{k=1}^{\infty} \gamma(k)=\infty.$  By definition   we obtain that
\begin{align*}
\beta(k)={\gamma(k) \over \gamma(k+1)}={(k+k_0+1)^{\delta}\over (k+k_0)^{\delta}}=\left( 1+{1\over k+k_0}\right)^{\delta}>1.
\end{align*}
Then   $\{\beta(k)\}$ is a monotonely decreasing sequence,  and $(1+{1/ k_0})^{\delta}=\beta(0) .$
Thus,  {\bf A6} (i) and (iii) hold.
\end{remark}

 Let $h \in \big(0, {2\over  \lambda_2(\mathbf{L})+\lambda_N(\mathbf{L})}\big)$ and   $\beta(0)  \in \big(1,  {1\over 1- h\lambda_2(\mathbf{L}})\big).$  We introduce some useful notations:
\begin{equation}\label{def3}
\begin{split}
 &\hat{\rho}_h \triangleq 1- h\lambda_2(\mathbf{L}), ~  \mathbf{z}_H \triangleq \left( z_i\mathbf{h}_1^T ,\dots, z_N\mathbf{h}_N^T \right)^T  ,\\ & M'(h,\beta(0)) \triangleq (1+2hd^*)\beta(0)+2h M_2(h,\beta(0)),  \mathcal{K}'( h,\beta(0))\triangleq \Big \lceil M'( h, \beta(0))-{1 \over 2} \Big \rceil
\end{split}
\end{equation}
with \begin{equation}\label{def-M1-M2}
\begin{split}
&M_1(h,\beta(0))\triangleq  \Big ( \sqrt{mN} C_x  (1+h \lambda_N(\mathbf{L}))   +
   {2 \|   \mathbf{z}_H\|_2  \over    \lambda_{\min}(\mathbf{F_d})}  \Big) \\& \times \Big(\| \mathbf{H_d}\|_{\infty}   +{
 h \lambda_N(\mathbf{L})\| \mathbf{H_d}\|_2  \over   1/    \beta(0) -\hat{\rho}_h  }\Big)  +  \| \mathbf{z}_H\|_{\infty} 
 + \lambda_N(\mathbf{L})
\Big(     \sqrt{mN} C_x\left(  1+h\beta(0)   \lambda_N(\mathbf{L})  \right)
   + { h   \| \mathbf{z}_H\|_2    \over  1 /  \beta(0)- \hat{\rho}_h}     \Big),
\\ &M_2(h,\beta(0))\triangleq \beta(0) \sqrt{mN} \lambda_N(\mathbf{L})   \Big( {    h  \lambda_N(\mathbf{L}) \over  2(1/  \beta(0) -    \hat{\rho}_h )}  
+{ 1  \over   \lambda_{\min}(\mathbf{F_d})}  \Big(\| \mathbf{H_d}\|_{\infty}   +{
h \lambda_N(\mathbf{L}) \| \mathbf{H_d}\|_2  \over   1/  \beta(0) -  \hat{\rho}_h  }\Big)  \Big),
\end{split}
\end{equation}
where $\mathbf{H_d}$ and $\mathbf{F_d}$ are defined in \eqref{def}.
   We now ready  to  state the main result of the algorithm \eqref{quantized-do}.

  \begin{proposition} \label{prp2} Suppose {\bf A4},  {\bf A5}, and {\bf A6} hold.
Let  Algorithm \ref{alg1}  be  applied to   the   least-squares  problem \eqref{LS}. Suppose  $h\in \left (0,\min \left \{{2\over \lambda_2(\mathbf{L})+\lambda_N(\mathbf{L})},{1\over \lambda_{\min}(\mathbf{F_d}) } \right\} \right)$ and $\beta(0)  \in \left(1,  {1\over 1- h\lambda_2(\mathbf{L})} \right).$    Then for any given $K \geq \mathcal{K}'( h,\beta(0))$,    the quantizer will never be saturated provided  that
\begin{align}
  s_r>  \max  \Big \{  & {C_x+h(C_x \| \mathbf{H_d}   \|_{\infty}  +\| \mathbf{z}_H\|_{\infty})   \over K+{1\over 2} },
 {M_1(h,\beta(0))/ M_2(h,\beta(0))} \Big \}. \label{ls-s0}
\end{align}  \end{proposition}

Proposition \ref{prp2} establishes the nonsaturation of the uniform quantizer, for which the proof is given in Section \ref{proof-thm02}.  Although the least-squares problem \eqref{LS}  seems like a special  case of  distributed optimization,  the main challenge lies in   that gradients of the quadratic function associated with each node  cannot be assumed to be globally bounded  a priori, a key technical assumption for the convergence analysis  of distributed (sub)gradient optimization \cite{nedic10,peng2014}. This is because the gradient function takes a linear form of the generated sequence $\mathbf{x}(k)$, which  might be unbounded with inappropriate algorithmic parameters.  Thus, the main  effort of the proof lies in  suitably choosing the  parameters and  proving  the boundedness of the generated sequence.

\begin{theorem}[High Data Rate]\label{thm3}
 Suppose {\bf A4},  {\bf A5}, and {\bf A6} hold.
Let    $h\in \left (0,\min \big \{{2\over \lambda_2(\mathbf{L})+\lambda_N(\mathbf{L})},{1\over \lambda_{\min}(\mathbf{F_d}) } \big\} \right)$ and $\beta(0)  \in \left(1,   1/( 1- h\lambda_2(\mathbf{L})) \right).$
Then for any given $K \geq \mathcal{K}'( h,\beta(0))$, along  Algorithm \ref{alg1}  there hold:
\begin{align}
 & \lim_{k\to \infty}  \mathbf{x}_{i}(k) =\mathbf{y_{\rm LS}}^*  \triangleq (\mathbf{H}\T\mathbf{H})^{-1}\mathbf{H}^T\mathbf{z}\quad \forall i \in \mathcal{V}, \label{ls-converg}
\\&\limsup_{k\rightarrow \infty}~ {   \| \mathbf{x}_i(k)-\mathbf{y_{\rm LS}}^*  \|_{\infty}\over \gamma(k)} <\infty\label{ls-rate}
\end{align}
 \end{theorem} provided that   $s_r$   satisfies  \eqref{ls-s0}.

 Theorem \ref{thm3}  shows that the Algorithm  \ref{alg1} can ensure  asymptotic convergence to the unique  least-squares solution $\mathbf{y_{\rm LS}}^*$.  Its  proof is deferred to Section  \ref{proof-thm3}.

\begin{remark} Note by Theorem \ref{thm3} that slow rate of convergence   is obtained by Algorithm \ref{alg1} with  decreasing  step-sizes for the  least-squares solver,   as opposed to  the  exponential   convergence of  the exact solver shown in Theorem \ref{thm1} for Algorithm \ref{alg1} with constant step-size. This is mainly because for the distributed least-squares problem even with un-quantized communication channel,  the primal domain algorithm cannot guarantee exact convergence with constant step-size \cite{shi2015extra}.  While it is noticed by  \cite{lei2016primal} and \cite{liu2017arrow} that the exact convergence or even the linear rate of convergence can  be obtained by  primal-dual domain algorithms.  As such, we might be able to  find the least-squares solution with limited
communication data rate at an exponential rate by the primal-dual domain methods.  We leave the problem
of  designing least-squares solver  with   non-decreasing step-size for future  research.
\end{remark}

Similar to Theorem \ref{thm2} for the exact solver case,  in the following theorem we show that
 we can also design a distributed protocol  for the   least-squares  solver to converge to a   least-squares  solution with $3-$level quantizers, which uses the minimum number of quantization levels.
   \begin{theorem}[Low Data Rate]\label{thm4}
 Suppose {\bf A4}, {\bf A5}  and  {\bf A6} hold.
Then  the following hold: \\
(i) For any $K\geq 1,$ $\Xi'_K$ is nonempty with
\begin{equation}\label{def_omega2}
\begin{split}
  \Xi'_K \triangleq& \Big  \{  ( h,\beta(0)):  \beta(0)  \in \big(1,   1 /(1- h\lambda_2(\mathbf{L})) \big) ,
\\&
h\in \Big (0,\min \big\{{2\over \lambda_2(\mathbf{L})+ \lambda_N(\mathbf{L})},{1\over \lambda_{\min}(\mathbf{F_d}) } \big \} \Big), 
 M'(h, \beta(0)) \leq  K+{1\over 2}\Big \}.
\end{split}
\end{equation} \\
(ii)   For any  $K\geq 1$, let $  ( h,\beta(0)) \in \Xi'_K$ and $s_r $  satisfy  \eqref{ls-s0}. Then  along  Algorithm \ref{alg1}  there hold for all $ i \in \mathcal{V} $ that
$ \lim_{k\to \infty}  \mathbf{x}_{i}(k) =\mathbf{y_{\rm LS}}^*   $ with the  rate  of convergence characterized by
$$\limsup_{k\rightarrow \infty}~ {   \| \mathbf{x}_i(k)-\mathbf{y_{\rm LS}}^*  \|_{\infty}\over \gamma(k)} <\infty.$$
\end{theorem}

The proof of Theorem \ref{thm4} is given Section  \ref{proof-thm4}.
 Similarly to Proposition ref{lem-rate},   the following result    with  the proof   given in Section \ref{proof-lem3}   gives an explicit  method  for   choosing  algorithm parameters  $( h,\beta(0)) \in \Xi_K'$ for any given $K \geq 1$ by  introducing a free parameter $\epsilon \in (0,1)$.
\begin{proposition} \label{lem-rate2}   For any given $K\geq 1 $  and   $\epsilon\in (0,1)$,  define $
h_{K,\epsilon}^*\triangleq  \min \left\{ {2\over \lambda_2(\mathbf{L})+ \lambda_N(\mathbf{L}},{1\over \lambda_{\min}(\mathbf{P_d}) } ,\hat{h}_{K,\epsilon} \right\}  $  and  $ \Xi'_{K,\epsilon} \triangleq   \Big \{ (\beta(0),h ):  \beta(0)^{-1}= 1-(1-\epsilon) h \lambda_{2}(\mathbf{L})  , h\in (0,h_{K,\epsilon}^*) \Big \}  ,  $
where $\hat{h}_{K,\epsilon}$ is defined in the following
\begin{align} \label{def-hath2}  & \hat{h}_{K,\epsilon}  \triangleq  2K \epsilon  \lambda_{\min}(\mathbf{F_d}) \Big( 2d^*\epsilon  \lambda_{\min}(\mathbf{F_d})  
+ (2K+1) \epsilon   (1-\epsilon)   \lambda_{\min}(\mathbf{F_d})\lambda_{2}(\mathbf{L})
 \notag  \\&
 +2\sqrt{mN} \lambda_N(\mathbf{L}) \times   \left(2 \epsilon   \| \mathbf{H_d}\|_{\infty}+  \kappa_N \big(2\| \mathbf{H_d}\|_2+\lambda_{\min}(\mathbf{F_d}) \big)\right) \Big)^{-1}
  \end{align}
with $\kappa_N  \triangleq { \lambda_N(\mathbf{L})  \over  \lambda_2(\mathbf{L}) }$.
Then $\Xi'_K=\bigcup_{\epsilon \in (0,1)} \Xi'_{K,\epsilon} $.
\end{proposition}

\subsection{Numerical Examples}\label{subsec:exact:ls:quantizaiton}


\noindent {\bf Example 4} Let $\mathbf{H},\mathbf{z}$ be given as follows:
\begin{equation*}
 \mathbf{H}=\left(
     \begin{array}{cc}
      1.7889 &	-1.0764\\
-1.0764	&0.1903\\
0.4707	&0.1008\\
0.8356	&-0.1716\\
0.5978	&-1.6668
     \end{array}
   \right),
   \mathbf{z}=\left(
       \begin{array}{c}
         -0.2854\\
1.2038\\
1.1032\\
0.7088\\
-0.9495
       \end{array}
     \right),
\end{equation*}
then the unique least square solution of $\mathbf{y}^*=\arg\min|| \mathbf{z}-\mathbf{H}\mathbf{y}||^2$ is
$\mathbf{y}^*=\left(
                    \begin{array}{c}
                      0.1415\\
                      0.6391
                    \end{array}
                  \right). $
The nodes again communicate according to the graph shown  in Fig \ref{fig_communication_graph}.

\begin{figure}[htbp]
\centering
\begin{minipage}[t]{0.45\textwidth}
\centering
  \includegraphics[width=3.2in]{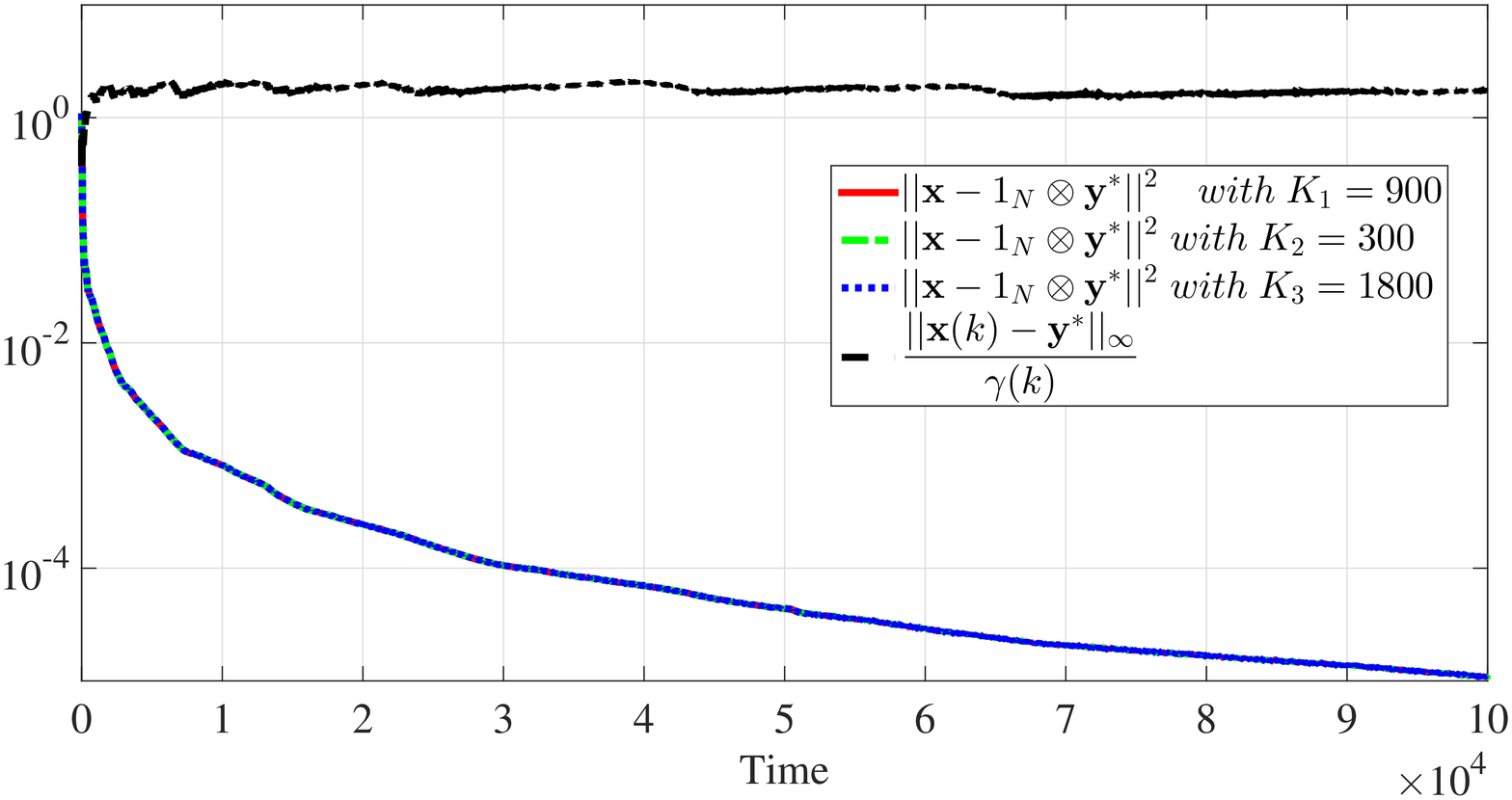}\\
  \caption{\small The trajectories of the sum of squared distance to the least square solution under $K=300, 900, 1800$.}\label{fig_ls_thm4_1}
\end{minipage}
\begin{minipage}[t]{0.08\textwidth}
\end{minipage}
\begin{minipage}[t]{0.45\textwidth}
\centering
  \includegraphics[width=3.2in]{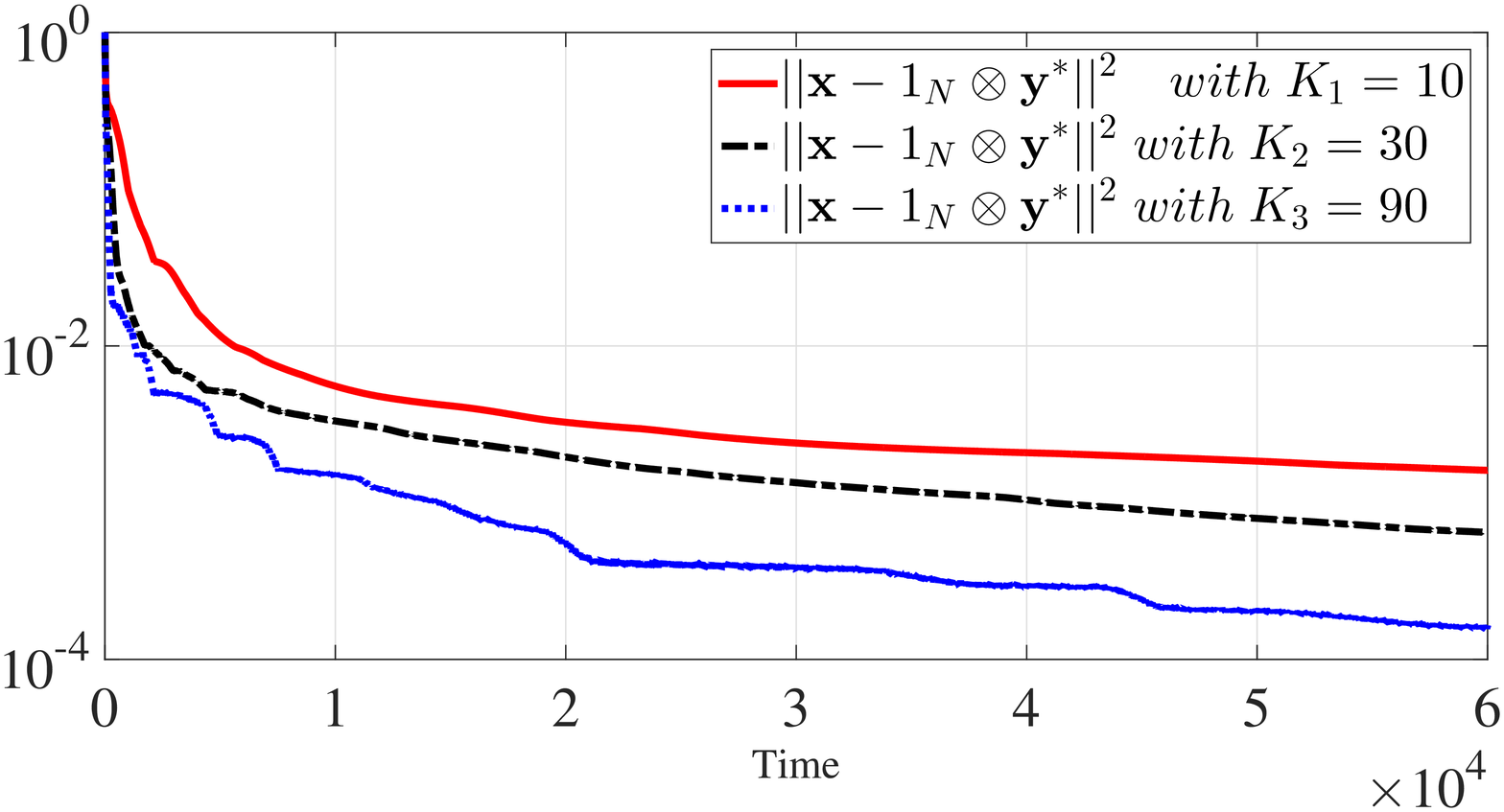}\\
  \caption{\small Trajectories of
$ ||\mathbf{x}-\mathbf{1}_N \otimes \mathbf{y}^* ||^2$  for  $K=20,50,100$ with the  algorithm parameters  chosen in Table \ref{table_1}. }\label{fig_ls_thm5_1}
\end{minipage}
\end{figure}

\noindent {\bf  [Validation of  Theorem \ref{thm3}.] }
Set $h=0.0853$  and   $\gamma(k)=(\frac{26}{k+26})^{0.85}$ such that $\beta(0)\in (1,\hat{\rho}_h^{-1})$.
Hence, $\mathcal{K}(h,\beta(0))=870$. We set $K_1=900,$ $K_2=300$ and $K_3=1800$, respectively.
We set $s_r=0.82$ to meet \eqref{ls-s0} fin all three cases.
We then run  Algorithm \ref{alg1} with the quantization levels $K_1$ $K_2$ and $K_3$, respectively,  while with the same parameters $h,\gamma(k)$. The simulation results are displayed in   Figure \ref{fig_ls_thm4_1}, which shows that the trajectories of $||\mathbf{x}-1_N\otimes \mathbf{y} ||^2$  coincide  in all  three cases. It then  implies that i) Once the sufficient condition of Theorem \ref{thm3} is satisfied, increasing data rate solely cannot speed up convergence;
ii) The condition in Theorem \ref{thm3} is sufficient for convergence but is not necessary.
Figure \ref{fig_ls_thm4_1} also shows the  trajectory of $\frac{||\mathbf{x}(k)-\mathbf{y}^* ||_{\infty}}{\gamma(k)}$, which verifies the convergence rate  described by  \eqref{ls-rate}.


\noindent {\bf  [Validation of  Theorem \ref{thm4}.] } We set the quantization level $K$ to be $K_1=10$, $K_2=30$ and $K_3=90$, respectively. Then we utilize Proposition  \ref{lem-rate2} to select algorithm parameters $h$ and $s(k)=s_r\gamma(k)=\frac{s_rk_0^{\delta}}{(k+k_0)^{\delta}}$ such that $(h,\beta(0))\in \Xi_{K}^{'}$ and $s_r$ satisfies \eqref{ls-s0} for the three cases. By setting $\epsilon=0.5$, the derived parameters for the three cases are given in Table \ref{table_1}.
Figure \ref{fig_ls_thm5_1} shows the trajectories of $||\mathbf{x}-\mathbf{1}_N\otimes \mathbf{y}^* ||^2$ for the three cases.
It demonstrates the convergence of the algorithm with the chosen parameters, verifying Theorem \ref{thm4}.
It also shows that with a higher data rate, the convergence could be faster if algorithm parameters are properly chosen.


\begin{table}[h!]
\centering
\begin{tabular}{|c|c|c|c|c|}
  \hline
  \quad     & $k_0$ & $\delta$ & $h$   & $s_r$ \\\hline
  K=10      & 120    & 0.85    & 0.0055 &  0.9583  \\\hline
  K=30      & 36     & 0.75    & 0.0164 &  0.6934 \\\hline
  K=90      & 9      & 0.55    & 0.0492 &  0.6968   \\\hline
\end{tabular}
\caption{Parameter settings}\label{table_1}
\end{table}
%

\subsection{Proofs of Statements}\label{Sec2:proof}

\subsubsection{ Preliminary Lemmas}
The following lemma gives a new but equivalent recursion  of Algorithm \ref{alg1}.
\begin{lemma} \label{ls-recur}  Let {\bf A4} and {\bf A6} (ii) hold.     Define
\begin{equation}\label{def2}
\begin{split} &\mathbf{P}  (k) =  \mathbf{I}_{mN}- h\left(  \mathbf{L}\otimes \mathbf{I}_m
+\gamma(k)\mathbf{H_d}\right)  ,\\&\bm{\varepsilon}(k) \triangleq {\mathbf{e}(k)/ s(k)},~  \bm{\eta}(k)\triangleq \left(\mathbf{D} \otimes \mathbf{I}_m \right)  { \mathbf{x}(k)  / \gamma(k)}  ,\end{split}
\end{equation}
 where $\mathbf{D} \triangleq \mathbf{I}_N -{\mathbf{1}_N\mathbf{1}_N ^T \over N}$ and $\mathbf{e}(k)$ is defined by \eqref{def}.    Then
\begin{align}
&\mathbf{x}(k+1) =   \mathbf{P}  (k)  \mathbf{x}(k) 
 +  h\gamma(k) \big(   s_r \mathbf{L}\otimes \mathbf{I}_m  \bm{\varepsilon}(k) + \mathbf{z}_H \big),  \label{ls-x}\\
&\bm{\eta}(k+1)  =  \beta(k) \big(  \left(\mathbf{I}_{mN}-  h  \mathbf{L}\otimes \mathbf{I}_m  \right) \bm{\eta}(k)
 +   hs_r   \mathbf{L} \otimes \mathbf{I}_m  \bm{\varepsilon}(k) +  h \mathbf{D} \otimes \mathbf{I}_m  \left( \mathbf{z}_H-\mathbf{H_d}\mathbf{x}(k)  \right) \big) \label{ls-eta},  \\
 &\bm{\varepsilon}(k+1)=\beta(k)  \big(  \bm{\theta}(k)- Q_K\left(  \bm{\theta}(k)\right)\big), \label{ls-ve}
\end{align}
where   $\bm{\theta}(k) $ is defined as follows:
\begin{equation}\label{def-ls-theta}
\begin{split}
 \bm{\theta}(k) &\triangleq  \big(\mathbf{I}_{mN} +h  \mathbf{L}\otimes \mathbf{I}_m   \big)\bm{\varepsilon}(k)
-    hs_r^{-1} \big(  \mathbf{L}\otimes \mathbf{I}_m   \bm{\eta}(k) + \mathbf{H_d}  \mathbf{x}(k) -   \mathbf{z}_H\big) . \end{split}
\end{equation}
 \end{lemma}
{\em Proof.}
By  using  \eqref{consensus} and   $\mathbf{e}(k)=s(k)\bm{\varepsilon}(k) =s_r\gamma(k)\bm{\varepsilon}(k)    $, we obtain  the following variant of \eqref{quantized-do}:
\begin{equation}\label{cf-x2}
\begin{split}
&\mathbf{x}(k+1)=\mathbf{x}(k)- h  \mathbf{L}\otimes \mathbf{I}_m  \mathbf{x}(k) + h  \mathbf{L}\otimes \mathbf{I}_m \mathbf{e}(k) 
- h\gamma(k)  \left( \mathbf{H_d} \mathbf{x}(k)-\mathbf{z}_H \right)
\\& =\left(\mathbf{I}_{mN}-  h  \mathbf{L}\otimes \mathbf{I}_m  \right)  \mathbf{x}(k) + hs_r\gamma(k)  \mathbf{L}\otimes \mathbf{I}_m \bm{\varepsilon}(k) 
 -h\gamma(k)   \mathbf{H_d} \mathbf{x}(k)+h\gamma(k) \mathbf{z}_H.
\end{split}
\end{equation}
Hence \eqref{ls-x} holds by using the definition of $ \mathbf{P}  (k) $ in \eqref{def2}.  By  multiplying both sides of   \eqref{cf-x2}  on  the left with ${\left(\mathbf{D}\otimes \mathbf{I}_m\right)     \over \gamma(k+1)} $,  using   $ \mathbf{D}  \mathbf{L}= \mathbf{L},\mathbf{D} (\mathbf{I}_N- h\mathbf{L})= (\mathbf{I}_N- h\mathbf{L}) \mathbf{D}$, and the definition  of  $\bm{\eta}(k)$ and $\beta(k) = {\gamma(k) \over \gamma(k+1)}  $, we obtain \eqref{ls-eta}.
By  subtracting  $\mathbf{ b}(k)$ from both sides of the first equality of   \eqref{cf-x2},
using $ \mathbf{e}(k)= \mathbf{x}(k)- \mathbf{b}(k)$, $ \gamma(k) =s_r^{-1} s(k)$   and $  \mathbf{L}\otimes \mathbf{I}_m  \mathbf{x}(k) =  \mathbf{L}\otimes \mathbf{I}_m  \left(\mathbf{D}\otimes \mathbf{I}_m\right)  \mathbf{x}(k) =\gamma(k)  \mathbf{L}\otimes \mathbf{I}_m \bm{\eta}(k)$, we obtain  that
 \begin{align*}
&{ \mathbf{x}(k+1)-\mathbf{b}(k)  } = \big(\mathbf{I}_{mN} +h \mathbf{L}\otimes \mathbf{I}_m  \big){ \mathbf{e}(k)  }
- h {  \mathbf{L}\otimes \mathbf{I}_m  \mathbf{x}(k) } - h{\gamma(k)  }(  \mathbf{H_d}  \mathbf{x}(k) - h  \mathbf{z}_H ) \\&= s(k) \big(\mathbf{I}_{mN} +h  \mathbf{L}\otimes \mathbf{I}_m   \big)\bm{\varepsilon}(k)
-    hs_r^{-1}s(k) \big(  \mathbf{L}\otimes \mathbf{I}_m   \bm{\eta}(k) + \mathbf{H_d}  \mathbf{x}(k) -   \mathbf{z}_H\big)   {{\tiny \eqref{def-ls-theta}} \atop =} s(k) \bm{\theta}(k) .
\end{align*}
  Recalling that  $\mathbf{ e}(k+1)=\mathbf{x}(k+1)- \mathbf{ b}(k+1)$ together with  \eqref{compact-form-xi},  we  then have the following:
\begin{equation*}
\begin{split}
&\mathbf{ e}(k+1)= \mathbf{x}(k+1)-\mathbf{ b}(k)  -s(k)  Q_K\left({\mathbf{x}(k+1)-\mathbf{ b}(k)\over s(k)} \right) =s(k)\big(  \bm{\theta}(k)- Q_K\left(  \bm{\theta}(k)\right)\big) ,
\end{split}
\end{equation*}
and hence dividing both sides of the above equation by $s(k+1)$  we obtain  \eqref{ls-ve}.
\hfill$\square$

\subsubsection{Proof of Proposition \ref{prp2}}\label{proof-thm02}
 The proof of  non-saturation of the uniform   quantizer is equivalent to show that for any $k\geq 0$, $ \bm{\theta}(k)$ defined  by \eqref{def-ls-theta} satisfies  $\| \bm{\theta}(k)\|_{\infty}<  K+{1\over 2}$. Again, we use an induction proof.

Recalling that $\gamma(0)=1, \mathbf{ b}_i(0)=0 ~ \forall i\in V,$  we obtain $\mathbf{ e}(0)= \mathbf{ x}(0)$ and $\bm{\varepsilon}(0)= \mathbf{ x}(0)/s_r.$  Then  by using $  \mathbf{L}\otimes \mathbf{I}_m \bm{\eta}(0)=   \mathbf{L}\otimes \mathbf{I}_m    \mathbf{x}(0)   $ and   \eqref{def-ls-theta}, we  obtain that
  \begin{align*}
&\bm{\theta}(0)=   \bm{\varepsilon}(0)+ h  \mathbf{L}\otimes \mathbf{I}_m    \bm{\varepsilon}(0)-    hs_r^{-1}   \mathbf{L}\otimes \mathbf{I}_m   \bm{\eta}(0)
 -hs_r^{-1} \big(  \mathbf{H_d}  \mathbf{x}(0) -   \mathbf{z}_H\big)  =  \mathbf{ x}(0)/s_r  - hs_r^{-1}(\mathbf{H_d}  \mathbf{x}(0)  -  \mathbf{z}_H)\big).
\end{align*}
   Then  from   {\bf A5},  and \eqref{ls-s0} it follows that
\begin{align*}
 &\| \bm{\theta}(0)\|_{\infty} \leq   C_x /s_r  + h\left( C_x \| \mathbf{H_d}   \|_{\infty}  +\| \mathbf{z}_H\|_{\infty}\right)/s_r
 = {C_x+h(C_x \| \mathbf{H_d}   \|_{\infty}  +\| \mathbf{z}_H\|_{\infty}) \over s_r}< K+{1\over 2}.
\end{align*}  Hence, when $k=0,$ the quantizer is unsaturated.  Next, for the induction, we assume  that when $k=0,\dots, p$, the quantizer is not saturated. Then by  \eqref{ls-ve}  and {\bf A6}, there holds
\begin{align}\label{ls-math-induct-z}
 \| \bm{\varepsilon}(k)  \|_{\infty} \leq {\beta(k) \over 2  } \leq  {\beta(0) \over 2 } \quad \forall k:1\leq k \leq p+1.
\end{align}

We aim to show that the quantizer is unsaturated for  $k=p+1$.
Define $\bm{\Gamma}(k ,k+1)\triangleq \mathbf{I}_{mN}$ and
\begin{align} \label{def-Gamma}
& \bm{\Gamma}(k_1,k_2)\triangleq   \mathbf{P} (k_1)  \mathbf{P} (k_1-1) \dots  \mathbf{P} (k_2)~\forall k_1\geq k_2\geq 0.
\end{align}
Then from \eqref{ls-x} and $\gamma(0)=1$ it follows that
\begin{equation}\label{ls-recursion-x}
\begin{split}
    \mathbf{x}(p+1)  &=\bm{\Gamma}(p,0) \mathbf{x}(0)+  h  s_r \bm{\Gamma}(p,1)   \mathbf{L}\otimes \mathbf{I}_m \bm{\varepsilon}(0)
  +  h s_r\sum_{i=1}^{p } \gamma(i)  \bm{\Gamma}(p,i+1)       \mathbf{L}\otimes \mathbf{I}_m \bm{\varepsilon}(i)  +h  \sum_{i=0}^{p }\gamma(i) \bm{\Gamma}(p,i+1)   \mathbf{z}_H.
\end{split}
\end{equation}

We  now  estimate the bound of $    \mathbf{x}(p+1).$ Using  $\mathbf{F_d}= \mathbf{L}\otimes \mathbf{I}_m + \mathbf{H_d}$,   the following holds:
\begin{align*}
& \min\{1, \gamma(k)\}x^T\mathbf{F_d}x \leq  x^T \left(  \mathbf{L}\otimes \mathbf{I}_m+\gamma(k)\mathbf{H_d} \right)x   \leq \max\{1, \gamma(k)\} x^T\mathbf{F_d}x   \quad \forall x\in \mathbb{R}^{mN}.
\end{align*}
By recalling that $0<\gamma(k) \leq 1$,  $\min\{1, \gamma(k)\}=\gamma(k)$ and $\max\{1, \gamma(k)\}=1.$ Thus,    for the matrix  $   \mathbf{L}\otimes \mathbf{I}_m + \gamma(k)\mathbf{H_d}  $, the smallest eigenvalue of is greater than or equal to $\gamma(k)\lambda_{\min}(\mathbf{F_d})$ while the   largest eigenvalue   is smaller than  or equal to $ \lambda_{\max}(\mathbf{F_d})$. Then by $ \mathbf{P}  (k) $ defined in \eqref{def2},
  the eigenvalues of $ \mathbf{P}  (k) $  sorted in an ascending order satisfy   $1- h \lambda_{\max}(\mathbf{F_d})\leq \lambda_1( \mathbf{P}  (k) )\leq \dots \leq \lambda_{mN}( \mathbf{P}  (k) )\leq 1-h \gamma(k)  \lambda_{\min}(\mathbf{F_d}).$  Thus, for any $k\geq 0:$
\begin{align*} \| \mathbf{P}  (k) \|_2&\leq \max \Big\{ \big |1-h \gamma(k)  \lambda_{\min}(\mathbf{F_d})\big|,  \big| 1- h  \lambda_{\max}(\mathbf{F_d})\big| \Big\} .
\end{align*}
Then by recalling that $0<h<{2\over \lambda_{\min}(\mathbf{F_d})+ \lambda_{\max}(\mathbf{F_d})}$ and $\gamma(k)\leq 1$,  the following holds:
 $$\| \mathbf{P}  (k) \|_2\leq  1-   h\gamma(k)\lambda_{\min}(\mathbf{F_d})  \leq \exp \left(-h\gamma(k)  \lambda_{\min}(\mathbf{F_d})\right) ,$$
where the last inequality holds by  $1-x\leq {\rm exp}(-x)~\forall x\geq 0.$
Then from \eqref{def-Gamma} it follows that for any $ k_1\geq k_2\geq 0:$
   $$\| \bm{\Gamma}(k_1,k_2)\|_2 <\exp \left(-h \lambda_{\min}(\mathbf{F_d}) \sum_{k=k_2}^{k_1}  \gamma(k) \right)
 .$$
Also, using  \eqref{ls-math-induct-z}, \eqref{ls-recursion-x}, {\bf A5}, $\|\mathbf{ x}\|_{\infty} \leq \|\mathbf{ x}\|_2\leq \sqrt{m} \|\mathbf{ x}\|_{\infty}$ for any $\mathbf{ x} \in \mathbb{R}^m$, and $ \mathbf{ x}(0)=s_r \bm{\varepsilon}(0),$  we obtain that
\begin{equation}\label{ls-recursion-x2}
\begin{split}
 &\| \mathbf{x}(p+1) \|_2    \leq    h s_r \|  \mathbf{L}\|_2\sum_{i=1}^{p } \gamma(i)   \| \bm{\Gamma}(p,i+1) \|_2    \|  \bm{\varepsilon}(i) \|_2
 +\| \bm{\Gamma}(p,0) \|_2 \| \mathbf{x}(0)\|_2\\&+  h  \| \bm{\Gamma}(p,1)\|_2  \|  \mathbf{L}\|_2 \| \mathbf{x}(0)\|_2  +h   \|  \mathbf{z}_H\|_2\sum_{i=1}^{p }\gamma(i) \| \bm{\Gamma}(p,i+1) \|_2
\\&\leq   \sqrt{mN} (1+h \lambda_N(\mathbf{L})) C_x
   +{hs_r \beta(0)\sqrt{mN}\lambda_N(\mathbf{L})  \over 2 }
\times \sum_{i=1}^{p }\gamma(i)   \exp \Big(- h\lambda_{\min}(\mathbf{F_d}) \sum_{k=i+1}^p\gamma(k)\Big)
   \\& + h \|   \mathbf{z}_H\|_2 \sum_{i=0}^{p }\gamma(i)\exp \big(-h \lambda_{\min}(\mathbf{F_d}) \sum_{k=i+1}^p\gamma(k)\big)  .
\end{split}
\end{equation}
Since  $h\gamma(k)  \leq  1/\lambda_{\min}(\mathbf{F_d})$ for any $k\geq 0$, there holds:
$$h\gamma(k) \leq 2\left(h \gamma(k)  -{\lambda_{\min}(\mathbf{F_d})\over 2}(h\gamma(k)) ^2\right) \quad \forall k\geq 0.$$
Thus, by $x-x^2/2<1-{\rm exp}(-x)~\forall x\in(0,1) $, we have the following sequence of inequalities:
\begin{align*}
&\sum_{i=k_1}^{p } h\lambda_{\min}(\mathbf{F_d})\gamma(i) \exp \Big(- h\lambda_{\min}(\mathbf{F_d}) \sum_{k=i+1}^p\gamma(k)\Big)
\\&\leq  2\sum_{i=k_1}^{p }   \Big( h\lambda_{\min}(\mathbf{F_d})\gamma(i)  -{\big(h\lambda_{\min}(\mathbf{F_d})\gamma(i)\big)^2/2} \Big) 
\times  \exp \Big(- h\lambda_{\min}(\mathbf{F_d}) \sum_{k=i+1}^p\gamma(i)\Big)
\\&\leq  2\sum_{i=k_1}^{p }\big( 1-  \exp \left(-\lambda_{\min}(\mathbf{F_d})\gamma(i)  \right)\big)
\times \exp \Big(-h \lambda_{\min}(\mathbf{F_d}) \sum_{k=i+1}^p\gamma(k)\Big)
 \\&= 2\sum_{i=k_1}^{p }\Big[ \exp \Big(- h\lambda_{\min}(\mathbf{F_d}) \sum_{k=i+1}^p\gamma(k)\Big)
-\exp \Big(-h \lambda_{\min}(\mathbf{F_d}) \sum_{k=i}^p\gamma(k)\Big) \Big]
\leq 2.
\end{align*}
 Using this in   \eqref{ls-recursion-x2} yields
\begin{equation}\label{ls-recursion-x3}
\begin{split}
 &  \| \mathbf{x}(p+1) \|_2 \leq    \sqrt{mN} C_x (1+h \lambda_N(\mathbf{L}))    +
 {s_r \beta(0) \sqrt{mN} \lambda_N(\mathbf{L})   \over     \lambda_{\min}(\mathbf{F_d})} + {2 \|   \mathbf{z}_H\|_2  \over    \lambda_{\min}(\mathbf{F_d})}\triangleq M_x.
\end{split}
\end{equation}

Since $\mathbf{L} $ is symmetric, we can define an orthogonal  matrix
$\mathbf{T}=\left({ \mathbf{1}_N \over \sqrt{N}}, \bm{\phi}_2,\dots,\bm{\phi}_N\right)$,
 where $  \mathbf{L}\bm{\phi}_i=\lambda_i(\mathbf{L}) \bm{\phi}_i$  for every $i=2,\dots,N$ .
Let $\tilde{\bm{\eta}}(k)=\left( \mathbf{T}^{-1} \otimes \mathbf{I}_m\right)\bm{\eta}(k)=\left(\mathbf{T}^T\otimes \mathbf{I}_m\right)\bm{\eta}(k)$  and decompose it as $\tilde{\bm{\eta}}(k)=(\tilde{\bm{\eta}}_1(k)^T,  \tilde{\bm{\eta}}_2(k)^T)^T$ with
$\tilde{\bm{\eta}}_1(k)={ 1 \over \sqrt{N}} \left(\mathbf{1}_N^T \otimes \mathbf{I}_m \right) \bm{\eta}  (k)$
and $\tilde{\bm{\eta}}_2(k)=\left( \mathbf{T}_2^T \otimes \mathbf{I}_m \right) \bm{\eta}  (k)$, where $ \mathbf{T}_2=\left( \bm{\phi}_2,\dots,\bm{\phi}_N\right)$.
Then $\tilde{\bm{\eta}}_1(k)=\mathbf{0}_m $ by $\bm{\eta}(k)= \left(\mathbf{D}\otimes \mathbf{I}_m\right)  { \mathbf{x}(k)   \over \gamma(k)}  $ and  $\mathbf{1}_N^T \mathbf{D}=\mathbf{0}_m^T.$
Then by multiplying both sides of  \eqref{ls-eta} with $ \mathbf{T}_2^T \otimes \mathbf{I}_m $ from the left,
 and noting  $ \mathbf{T}_2 ^T\mathbf{L}= {\rm diag} \left\{ \lambda_2(\mathbf{L}),\dots,  \lambda_N(\mathbf{L})\right \}\mathbf{T}_2 ^T $ we have the following:
\begin{align}
\tilde{\bm{\eta}}_2(k+1)=& \beta(k)   hs_r \left(\mathbf{T}_2^T\mathbf{L} \otimes \mathbf{I}_m\right) \bm{\varepsilon}(k)  + \beta(k)  h \left(\mathbf{T}_2^T\mathbf{D} \otimes \mathbf{I}_m\right) \left( \mathbf{z}_H-\mathbf{H_d}\mathbf{x}(k) \right)\notag
\\&+\beta(k) \Big( \underbrace{ {\rm diag} \left\{ 1- h\lambda_2(\mathbf{L}),\dots, 1- h\lambda_N(\mathbf{L})\right \}}_{\mathbf{D}_h} \otimes \mathbf{I}_m\Big)\tilde{\bm{\eta}}_2(k)  \label{def-Dh}
  \end{align}
  Thus, there holds
\begin{equation}\label{def-tilde-eta}
\begin{split}
\tilde{\bm{\eta}}_2(p+1)&= \left(\mathbf{D}_h^{p+1} \otimes \mathbf{I}_m\right)\tilde{\bm{\eta}}_2(0)\prod_{k=0}^{p}\beta(k)  +     hs_r \sum_{k=0}^p  \left(\mathbf{D}_h^{k}\mathbf{T}_2^T\mathbf{L} \otimes \mathbf{I}_m\right)  \prod_{i=p-k}^{p}\beta(i)\bm{\varepsilon}(p-k)
  \\&+ h\sum_{k=0}^p  \left(\mathbf{D}_h^{k}\mathbf{T}_2^T\mathbf{D}\otimes \mathbf{I}_m\right)  \prod_{i=p-k}^{p}\beta(i)\left( \mathbf{z}_H-\mathbf{H_d}\mathbf{x}(p-k) \right).  \end{split}
\end{equation}
 Note that $ \left( \mathbf{T}_2  \otimes \mathbf{I}_m \right) \tilde{\bm{\eta}}_2  (k)
 =\left( \mathbf{T}_2 \mathbf{T}_2 ^T \otimes \mathbf{I}_m \right) \bm{\eta}   (k) =
 \left( \mathbf{I}_N-{\mathbf{1}_N\mathbf{1}_N^T \over N} \otimes \mathbf{I}_m \right) \bm{\eta}   (k)= \bm{\eta}   (k).$
 Then by multiplying both sides of  \eqref{def-tilde-eta}  on the left with $\left( \mathbf{T}_2  \otimes \mathbf{I}_m \right) $, there holds
\begin{equation}\label{def-etap}
\begin{split}
 \bm{\eta} (p+1) &= \left(\mathbf{T}_2\mathbf{D}_h^{p+1} \mathbf{T}_2^T \otimes \mathbf{I}_m\right) \bm{\eta} (0)\prod_{k=0}^{p}\beta(k)
 +     h s_r\sum_{k=0}^p \left( \mathbf{T}_2\mathbf{D}_h^{k}\mathbf{T}_2^T\mathbf{L} \otimes \mathbf{I}_m\right) \bm{\varepsilon}(p-k) \prod_{i=p-k}^{p}\beta(i)
  \\&+ h \sum_{k=0}^p \left( \mathbf{T}_2\mathbf{D}_h^{k}\mathbf{T}_2^T\mathbf{D}\otimes \mathbf{I}_m\right)
   \left( \mathbf{z}_H-\mathbf{H_d}\mathbf{x}(p-k) \right)\prod_{i=p-k}^{p}\beta(i) .  \end{split} \end{equation}
 By the definition of $\mathbf{D}_h$ in \eqref{def-Dh},
$\| \mathbf{D}_h\|_2=\max\{ | 1- h\lambda_2(\mathbf{L})|, |1- h\lambda_N(\mathbf{L}) |\}$. Thus, by using    $h \in (0, {2\over  \lambda_2(\mathbf{L})+\lambda_N(\mathbf{L})})$  and
    \cite[Lemma 3.1]{litao2011},  we obtain that $\| \mathbf{D}_h\|_2= 1- h\lambda_2(\mathbf{L})=\hat{\rho}_h.$
Taking  two-norms of \eqref{def-etap}, by recalling that  $\beta(k) \leq \beta(0) ~\forall k\geq 0$,
$ \|\mathbf{D}\|_2=\|\mathbf{T_2}\|_2 =1,$  we have  the following:
  \begin{align*}
 &\| \bm{\eta} (p+1)\|_2 \leq { \left(\beta(0)  \hat{\rho}_h \right)} ^{p+1} \| \bm{\eta} (0)\|_2 + hs_r  \beta(0)  \|\mathbf{L}\|_2 { \left(\beta(0)  \hat{\rho}_h \right)} ^{p} \| \bm{\varepsilon}(0)\|_2
  \\&+     hs_r  \beta(0) \|\mathbf{L}\|_2 \sum_{k=0}^{p-1} { \left(\beta(0)  \hat{\rho}_h \right)} ^{k} \| \bm{\varepsilon}(p-k)\|_2  +h \beta(0)\sum_{k=0}^p { \left(\beta(0)  \hat{\rho}_h \right)} ^{k} \left \| \mathbf{z}_H-\mathbf{H_d}\mathbf{x}(p-k) \right \|_2.  \end{align*}
Note  by  $\gamma(0)=1, s(0)=s_r$, $ \|\mathbf{D}\|_2=1$,   and {\bf A5}  that
\begin{align*}
& \| \bm{\eta} (0)\|_2 = \| \mathbf{D} \otimes \mathbf{I}_m \mathbf{x}(0)\|_2  \leq \| \mathbf{D}\|_2 \|\mathbf{x}(0)\|_2   \leq \sqrt{mN}C_x ,\\& \| \bm{\varepsilon}(0)\|_2\leq \|\mathbf{x}(0)\|_2/s(0) \leq
  \sqrt{mN}C_x/s_r  .
  \end{align*}
Similar to   \eqref{ls-recursion-x3} we can easily show that $\|\mathbf{x}(k)\|_2 \leq M_x ~\forall k=0,\dots,p.$
 Then by using    \eqref{ls-math-induct-z}, $\beta(0)  \hat{\rho}_h<1$ and $\sum_{k=0}^p { \left(\beta(0)  \hat{\rho}_h \right)} ^{k} \leq{1   \over  1- {  \beta(0)  \hat{\rho}_h  } } $, we obtain the following:
     \begin{align}\label{def-ls-meta}
  \| \bm{\eta} (p+1)\|_2  &\leq   \sqrt{mN} C_x\left(  1+h\beta(0)   \lambda_N(\mathbf{L})  \right) 
  +   {  \sqrt{mN}hs_r  \beta(0)^2   \lambda _N(\mathbf{L})  \over  2 (1 -\beta(0)\hat{\rho}_h )  }  \notag \\&
  + \left ( \| \mathbf{z}_H\|_2 + \|\mathbf{H_d}\|_2 M_x\right){ h \beta(0)   \over  1  -    \beta(0) \hat{\rho}_h}   .
  \end{align}
  This together with  \eqref{def-ls-theta},  \eqref{ls-math-induct-z} and  \eqref{ls-recursion-x3} leads to
\begin{align*}
&\| \bm{\theta}(p+1)\|_{\infty}  \leq  \|  \left( \mathbf{I}_{mN}+ h  \mathbf{L}\otimes \mathbf{I}_m \right) \bm{\varepsilon}(p+1) \|_{\infty} +
 hs_r^{-1}   \times \big (\lambda_N(\mathbf{L}) \|  \bm{\eta}(p+1)\|_2 +  \| \mathbf{z}_H\|_{\infty}+    \| \mathbf{H_d}\|_{\infty} \|  \mathbf{x}(p+1)\|_{\infty} \big)
 \\& \leq   \beta(0)(1/2+ hd^*)  + h \sqrt{mN} \lambda_N(\mathbf{L})  \Big (   \beta(0)      \lambda^{-1}_{\min}(\mathbf{F_d})    \times \Big(\| \mathbf{H_d}\|_{\infty}   +{
h\beta(0) \lambda_N(\mathbf{L}) \| \mathbf{H_d}\|_2  \over   1 -    \beta(0)\hat{\rho}_h  }\Big)  +{    h\beta(0)^2 \lambda_N(\mathbf{L}) \over  2 (1 -    \beta(0)\hat{\rho}_h )} \Big )
\\&  + hs_r^{-1}\Big ( \sqrt{mN} C_x  (1+h \lambda_N(\mathbf{L}))   +
   {2 \|   \mathbf{z}_H\|_2  \over    \lambda_{\min}(\mathbf{F_d})}  \Big)
 \times \Big(\| \mathbf{H_d}\|_{\infty}   +{
 h\beta(0) \lambda_N(\mathbf{L})\| \mathbf{H_d}\|_2  \over   1 -    \beta(0)\hat{\rho}_h  }\Big) +  hs_r^{-1}  \| \mathbf{z}_H\|_{\infty}     \\&  + hs_r^{-1} \lambda_N(\mathbf{L})
\Big(     \sqrt{mN} C_x\left(  1+h\beta(0)   \lambda_N(\mathbf{L})  \right)
+ { h \beta(0) \| \mathbf{z}_H\|_2    \over  1  -    \beta(0) \hat{\rho}_h}     \Big)
 \\&= { \beta(0)(1+2hd^*)\over 2  } + h\left(s_r^{-1} M_1(h,\beta(0))+ M_2(h,\beta(0)) \right) \qquad \scriptstyle{({\rm by~ }  \eqref{def-M1-M2}) }
\\& \leq  \beta(0)(1/2+hd^*) + 2hM_2(h,\beta(0))  = M'( h, \beta(0))~\scriptstyle{({\rm by ~}  \eqref{def3}{\rm ~and ~} \eqref{ls-s0})} \notag
\\& \leq \Big \lceil M'( h, \beta(0))-{1 \over 2} \Big \rceil +{1\over 2}=\mathcal{K}'( h,\beta(0))+{1 \over 2} \leq K+{1\over 2}.  \notag
\end{align*}
As a result, when $k=p+1$, the quantizer is also unsaturated.
Therefore, by induction, we conclude that if a $(2K+1)$-levels uniform quantizer is applied, then the  quantizer will   never be saturated.
\hfill $\blacksquare$
 \subsubsection{Proof of Theorem \ref{thm3}}\label{proof-thm3}
 From  Proposition \ref{prp2}  it follows  that  \eqref{def-ls-meta}  holds for any $p\geq 0.$
This   implies that $\sup_{k\geq 1} \|\bm{\eta}(k)\|_{\infty} <\infty.$ Then using  the  definition of  $\bm{\eta}(k)$, we obtain that $  \| \mathbf{x}_i(k)-\mathbf{y}(k) \|_{\infty} =\mathcal{O}( \gamma(k)).$

 Define $ \mathbf{y} _k=   \sum_{i=1}^N \mathbf{x}_{i,k}/N.$
 Then by multiplying both sides of \eqref{cf-x2} from the left by ${1\over N}  \left(\mathbf{1}_N\otimes \mathbf{I}_m\right) ,$
there holds
 \begin{align*}
  \mathbf{y} _{k+1}
& = \mathbf{y} _k-   h\gamma(k)   (\mathbf{H}\T\mathbf{H}  \mathbf{y} _k-\mathbf{H}^T\mathbf{z} )/N
 - {h \gamma(k) \over N} \sum_{i=1}^N   \mathbf{h}_i\mathbf{h}_i^\top  \left( \mathbf{x}_{i,k}-  \mathbf{y} _k\right).
 \end{align*}
 Then by recalling that $\mathbf{y}^\star=(\mathbf{H}\T\mathbf{H})^{-1}\mathbf{H}^T\mathbf{z}$, we obtain  that
 \begin{align*}
  \mathbf{y} _{k+1}- \mathbf{y_{\rm LS}}^*&= \mathbf{y} _k- \mathbf{y_{\rm LS}}^*
-    {h\gamma(k) \over N}  \mathbf{H}\T\mathbf{H} \left( \mathbf{y} _k- \mathbf{y_{\rm LS}}^*\right)
 - {h\gamma(k) \over N} \sum_{i=1}^N   \mathbf{h}_i\mathbf{h}_i^\top  \left( \mathbf{x}_{i,k}-  \mathbf{y} _k\right).
 \end{align*}
Since $ \left( \mathbf{x}_{i,k}-  \mathbf{y} _k\right) \to \mathbf{0}$ and $ \mathbf{H}\T\mathbf{H}$ is positive definite by ${\rm rank}( \mathbf{H})=m$,  from \cite[Lemma 3.1.1]{Chen_2002}  it follows  that $\lim\limits_{k\to \infty}\mathbf{y} _k= \mathbf{y_{\rm LS}}^*$, this together with $ \mathbf{x}_i(k)-\mathbf{y}(k) \to
 \mathbf{0}$ implies \eqref{ls-converg}.
Then by   $  \| \mathbf{x}_i(k)-\mathbf{y}(k) \|_{\infty} =\mathcal{O}( \gamma(k))$ and  \cite[Theorem 3.1.1]{Chen_2002} we obtain that
 $  \|  \mathbf{y}(k) -\mathbf{y_{\rm LS}}^* \|_{\infty} =\mathcal{O}( \gamma(k))$, and hence \eqref{ls-rate} holds.
 \hfill $\blacksquare$

\subsubsection{Proof of Theorem \ref{thm4}  }\label{proof-thm4}
\noindent (i) Using $\hat{\rho}_h=1- h\lambda_2(\mathbf{L}) $, $\kappa_N ={ \lambda_N(\mathbf{L})  \over  \lambda_2(\mathbf{L}) }$,  \eqref{def3} and \eqref{def-M1-M2},  we obtain the following:
\begin{equation*}
\begin{split}
&M'(h,1)
 = { 1+2hd^*\over 2 }+ 2h     \sqrt{mN} \lambda_N(\mathbf{L})  \Big(    \lambda^{-1}_{\min}(\mathbf{F_d})  \big(\| \mathbf{H_d}\|_{\infty}    +
\| \mathbf{H_d}\|_2 \kappa_N  \big)  +{  \kappa_N   \over  2  } \Big)
\end{split}
\end{equation*}
Thus,  $ \lim\limits_{h\to 0 }  M'(h ,1) ={1\over 2} $.
 Then for any $K\geq 1,$ there exists $h^*\in \left (0,\min \left \{{2\over \lambda_2(\mathbf{L})+\lambda_N(\mathbf{L})},{1\over \lambda_{\min}(\mathbf{P_d}) } \right\} \right)$ such that  $ M'(h^*, 1) \leq K$. By recalling  the definition of $ M'(h, \beta(0))$ in   \eqref{def3}, we have that
$$\lim_{\beta(0) \to 1} M'(h^* ,\beta(0)) =M'(h^*, 1) \leq K.$$
Then  there  exists  $ \beta^*(0)\in \left(1, {1\over 1- h\lambda_2(\mathbf{L}}\right) $ such that   $ M'(h^*, \beta^*(0)) \leq K+{1\over 2}$.
Therefore,   $(h^*,\beta^*(0))  \in \Xi'_K  $. Hence $\Xi'_K $ is nonempty.

\noindent (ii) For any $( h,\beta(0)) \in \Xi'_K$,  from  \eqref{def_omega2} it follows that $
h\in \left (0,\min \big \{{2\over \lambda_2(\mathbf{L})+\lambda_N(\mathbf{L})},{1\over \lambda_{\min}(\mathbf{P_d}) } \big\} \right),$  $ \beta(0)\in \left(1, {1\over 1- h\lambda_2(\mathbf{L}} \right),$ and $ M'(h,\beta(0)) \leq   K+{1\over 2}$. Then by definition  \eqref{def3},
$$  \mathcal{K}'( h,\beta(0))=\Big \lceil M'(h,\beta(0))-{1\over 2} \Big \rceil \leq K,$$ which  together with  Theorem \ref{thm3} leads to assertion (ii).
\hfill $\blacksquare$

 \subsubsection{  Proof of  Proposition  \ref{lem-rate2}  } \label{proof-lem3}
  We first validate $\bigcup_{\epsilon \in (0,1)} \Xi'_{K,\epsilon}\subset \Xi_K$.
For any given $K\geq 1 $ and $ \epsilon \in (0,1)$, let $(\beta(0),h)\in  \Xi'_{K,\epsilon}$.
Then  $1/\beta(0)-\hat{\rho}_h= \epsilon h \lambda_{2}(\mathbf{L})  >0$ by   $\hat{\rho}_h=1-  h \lambda_{2}(\mathbf{L})$, and  hence $\beta(0)<1/\hat{\rho}_h$. Then by  the definition of $  M'(h,\beta(0))$ in \eqref{def3}, using $\kappa_N ={ \lambda_N(\mathbf{L})  \over  \lambda_2(\mathbf{L}) }$ and $ \beta(0)^{-1}= 1-(1-\epsilon) h \lambda_{2}(\mathbf{L}) $, we obtain  that \begin{align*}
 &  M'(h,\beta(0))  =   (1+2hd^*)\beta(0) + 2h\sqrt{mN} \lambda_N(\mathbf{L}) \beta(0)  \times \Big( { 1  \over   \lambda_{\min}(\mathbf{F_d})}
  \Big(\| \mathbf{H_d}\|_{\infty}   +{  \kappa_N \| \mathbf{H_d}\|_2  \over   \epsilon}\Big)  +{    \kappa_N \over  2 \epsilon}\Big)    \\& = \left(2 \epsilon   \| \mathbf{H_d}\|_{\infty}+ 2\kappa_N \| \mathbf{H_d}\|_2
  +\kappa_N \lambda_{\min}(\mathbf{F_d})  \right)
  \times{ (1+2hd^*)\epsilon  \lambda_{\min}(\mathbf{F_d})  +2h\sqrt{mN} \lambda_N(\mathbf{L})
 \over 2 \epsilon (1-  (1-\epsilon)h \lambda_{2}(\mathbf{L}) )   \lambda_{\min}(\mathbf{F_d})}
\end{align*}
 Then by using the definition of $ \hat{h}_{K,\epsilon} $  in \eqref{def-hath2} and $h<\hat{h}_{K,\epsilon}$,  there holds  $   M'(h,\beta(0)) < K+{1\over 2}$.  Obviously, for any $(\alpha,h)\in  \Xi_{K,\epsilon}$,  there holds $ h\in \left (0,\min \big\{{2\over \lambda_2(\mathbf{L})+ \lambda_N(\mathbf{L})},{1\over \lambda_{\min}(\mathbf{F_d}) } \big \} \right)$.  In summary, we have  verified $\Xi'_{K,\epsilon}\subset \Xi'_K $   for any given $K\geq 1 $ and  any $ \epsilon \in (0,1)$. Thus, $\bigcup_{\epsilon \in (0,1)} \Xi'_{K,\epsilon}\subset \Xi'_K$.

 Similar to that of  Lemma \ref{lem-rate}, we can also  prove that  $\Xi'_K\subset \bigcup_{\epsilon \in (0,1)} \Xi'_{K,\epsilon}$. Thus, we complete  the proof.  \hfill $\blacksquare$

\section{Conclusions}\label{sec:conclusions}

We have studied solving  linear equations over a network subject to digital node communications with a limited data rate.
We propose  a  node encoder-decoder pair,  based on which a distributed quantized algorithm is designed. For the unique exact solution case, the proposed encoder-decoder powered algorithm drove  each node state to the solution asymptotically at an exponential rate. For the unique least-squares solution case, the same encoder-decoder pair enabled  the algorithm to compute  such a solution with a  properly  selected time-varying step size. A minimal data  rate was shown to be enough for  the desired convergence for both cases. These results suggest the practical applicability  of various network linear equation solvers  in the literature.


\begin{thebibliography}{10}

\bibitem{tsi}
J.~Tsitsiklis, D.~Bertsekas, and M.~Athans, ``Distributed asynchronous
  deterministic and stochastic gradient optimization algorithms,'' \emph{IEEE
  transactions on automatic control}, vol.~31, no.~9, pp. 803--812, 1986.

\bibitem{jad03}
A.~Jadbabaie, J.~Lin, and A.~S. Morse, ``Coordination of groups of mobile
  autonomous agents using nearest neighbor rules,'' \emph{IEEE Transactions on
  automatic control}, vol.~48, no.~6, pp. 988--1001, 2003.

\bibitem{xiao04}
L.~Xiao and S.~Boyd, ``Fast linear iterations for distributed averaging,''
  \emph{Systems \& Control Letters}, vol.~53, no.~1, pp. 65--78, 2004.

\bibitem{Rabbat2010}
A.~G. Dimakis, S.~Kar, J.~M. Moura, M.~G. Rabbat, and A.~Scaglione, ``Gossip
  algorithms for distributed signal processing,'' \emph{Proceedings of the
  IEEE}, vol.~98, no.~11, pp. 1847--1864, 2010.

\bibitem{magnusbook}
M.~Mesbahi and M.~Egerstedt, \emph{Graph theoretic methods in multiagent
  networks}.\hskip 1em plus 0.5em minus 0.4em\relax Princeton University Press,
  2010.

\bibitem{mou13}
S.~Mou and A.~Morse, ``A fixed-neighbor, distributed algorithm for solving a
  linear algebraic equation,'' in \emph{European Control Conference (ECC), } pp. 2269--2273, 2013.

\bibitem{rabbat2004}
M.~Rabbat and R.~Nowak, ``Distributed optimization in sensor networks,'' in
  \emph{Proceedings of the 3rd international symposium on Information
  processing in sensor networks}.\hskip 1em plus 0.5em minus 0.4em\relax ACM,
  2004, pp. 20--27.

\bibitem{nedic09}
A.~Nedic and A.~Ozdaglar, ``Distributed subgradient methods for multi-agent
  optimization,'' \emph{IEEE Transactions on Automatic Control}, vol.~54,
  no.~1, pp. 48--61, 2009.

\bibitem{yi2016initialization}
P.~Yi, Y.~Hong, and F.~Liu, ``Initialization-free distributed algorithms for
  optimal resource allocation with feasibility constraints and application to
  economic dispatch of power systems,'' \emph{Automatica}, vol.~74, pp.
  259--269, 2016.

\bibitem{margaris2014parallel}
A.~Margaris, S.~Souravlas, and M.~Roumeliotis, ``Parallel implementations of
  the jacobi linear algebraic systems solve,'' \emph{arXiv:1403.5805}, 2014.

\bibitem{saad1999distributed}
Y.~Saad and M.~Sosonkina, ``Distributed schur complement techniques for general
  sparse linear systems,'' \emph{SIAM Journal on Scientific Computing},
  vol.~21, no.~4, pp. 1337--1356, 1999.

\bibitem{anderson1997}
C.~Andersson, ``Solving linear eqauations on parallel distributed memory
  architectures by extrapolation,'' \emph{Technical Report, Royal Institute of
  Technology}, 1997.

\bibitem{mehmood05}
R.~Mehmood and J.~Crowcroft, ``Parallel iterative solution method for large
  sparse linear equation systems,'' University of Cambridge, Computer
  Laboratory, Tech. Rep., 2005.

\bibitem{lei2015distributed}
J.~Lei and H.-F. Chen, ``Distributed randomized pagerank algorithm based on
  stochastic approximation,'' \emph{IEEE Transactions on Automatic Control},
  vol.~60, no.~6, pp. 1641--1646, 2015.

\bibitem{nedic10}
A.~Nedic, A.~Ozdaglar, and P.~A. Parrilo, ``Constrained consensus and
  optimization in multi-agent networks,'' \emph{IEEE Transactions on Automatic
  Control}, vol.~55, no.~4, pp. 922--938, 2010.

\bibitem{elia}
J.~Wang and N.~Elia, ``Control approach to distributed optimization,'' in
  \emph{The 48th Annual Allerton Conference on Communication, Control, and Computing (Allerton),}  pp. 557--561, 2010.

\bibitem{lei2016primal}
J.~Lei, H.-F. Chen, and H.-T. Fang, ``Primal--dual algorithm for distributed
  constrained optimization,'' \emph{Systems \& Control Letters}, vol.~96, pp.
  110--117, 2016.

\bibitem{lu12}
J.~Lu and C.~Y. Tang, ``Zero-gradient-sum algorithms for distributed convex
  optimization: The continuous-time case,'' \emph{IEEE Transactions on
  Automatic Control}, vol.~57, no.~9, pp. 2348--2354, 2012.

\bibitem{Mou-TAC-2015}
S.~Mou, J.~Liu, and A.~S. Morse, ``A distributed algorithm for solving a linear
  algebraic equation,'' \emph{IEEE Transactions on Automatic Control}, vol.~60,
  no.~11, pp. 2863--2878, 2015.

\bibitem{Shi-TAC-LAE}
G.~Shi, B.~D. Anderson, and U.~Helmke, ``Network flows that solve linear
  equations,'' \emph{IEEE Transactions on Automatic Control}, vol.~62, no.~6,
  pp. 2659--2674, 2017.

\bibitem{Yang-IFAC-2017}
Y.~Liu, C.~Lageman, B.~D. Anderson, and G.~Shi, ``Exponential least squares
  solvers for linear equations over networks,'' \emph{IFAC World Congress},
  vol.~50, no.~1, pp. 2543--2548, 2017.

\bibitem{Yang-CDC-2017}
Y.~Liu, Y.~Lou, B.~D.~O. Anderson, and G.~Shi, ``Network flows as least squares
  solvers for linear equations,'' in \emph{IEEE Conference on Decision and Control,}  pp. 1046--1051,  2017.

\bibitem{brian15}
B.~Anderson, S.~Mou, A.~S. Morse, and U.~Helmke, ``Decentralized gradient
  algorithm for solution of a linear equation,'' \emph{Numerical Algebra, Control \& Optimization},  vol.~6, no.~3,
 pp. 319--328, 2016.

\bibitem{jadbabaie15}
R.~Tutunov, H.~B. Ammar, and A.~Jadbabaie, ``A fast distributed solver for
  symmetric diagonally dominant linear equations,'' \emph{arXiv:1502.03158}, 2015.

\bibitem{asuman14}
C.~E. Lee, A.~Ozdaglar, and D.~Shah, ``Solving systems of linear equations:
  Locally and asynchronously,'' \emph{Computing Research Repository},   2014.

\bibitem{jvn49}
J.~Von~Neumann, ``On rings of operators. reduction theory,'' \emph{Annals of
  Mathematics}, pp. 401--485, 1949.

\bibitem{shitac}
G.~Shi, K.~H. Johansson, and Y.~Hong, ``Reaching an optimal consensus:
  Dynamical systems that compute intersections of convex sets,'' \emph{IEEE
  Transactions on Automatic Control}, vol.~58, no.~3, pp. 610--622, 2013.

\bibitem{Brockett2000}
R.~W. Brockett and D.~Liberzon, ``Quantized feedback stabilization of linear
  systems,'' \emph{IEEE transactions on Automatic Control}, vol.~45, no.~7, pp.
  1279--1289, 2000.

\bibitem{basar}
A.~Kashyap, T.~Ba{\c{s}}ar, and R.~Srikant, ``Quantized consensus,''
  \emph{Automatica}, vol.~43, no.~7, pp. 1192--1203, 2007.

\bibitem{carli2009}
P.~Frasca, R.~Carli, F.~Fagnani, and S.~Zampieri, ``Average consensus on
  networks with quantized communication,'' \emph{International Journal of
  Robust and Nonlinear Control}, vol.~19, no.~16, pp. 1787--1816, 2009.

\bibitem{girish2011}
G.~N. Nair, F.~Fagnani, S.~Zampieri, and R.~J. Evans, ``Feedback control under
  data rate constraints: An overview,'' \emph{Proceedings of the IEEE},
  vol.~95, no.~1, pp. 108--137, 2007.

\bibitem{hong2016}
Z.~Qiu, L.~Xie, and Y.~Hong, ``Quantized leaderless and leader-following
  consensus of high-order multi-agent systems with limited data rate,''
  \emph{IEEE Transactions on Automatic Control}, vol.~61, no.~9, pp.
  2432--2447, 2016.

\bibitem{litao2011}
T.~Li, M.~Fu, L.~Xie, and J.-F. Zhang, ``Distributed consensus with limited
  communication data rate,'' \emph{IEEE Transactions on Automatic Control},
  vol.~56, no.~2, pp. 279--292, 2011.

\bibitem{rabbat2005}
M.~G. Rabbat and R.~D. Nowak, ``Quantized incremental algorithms for
  distributed optimization,'' \emph{IEEE Journal on Selected Areas in
  Communications}, vol.~23, no.~4, pp. 798--808, 2005.

\bibitem{peng2014}
P.~Yi and Y.~Hong, ``Quantized subgradient algorithm and data-rate analysis for
  distributed optimization,'' \emph{IEEE Transactions on Control of Network
  Systems}, vol.~1, no.~4, pp. 380--392, 2014.

\bibitem{lei2018datatare}
J.~Lei, P.~Yi, G.~Shi, and B.~D.~O. Anderson, ``Network linear equations with
  finite data rates,'' \emph{the Proceedings of the IEEE Conference on Decision
  and Control}, 2018.

\bibitem{shi2015extra}
W.~Shi, Q.~Ling, G.~Wu, and W.~Yin, ``Extra: An exact first-order algorithm for
  decentralized consensus optimization,'' \emph{SIAM Journal on Optimization},
  vol.~25, no.~2, pp. 944--966, 2015.

\bibitem{liu2017arrow}
Y.~Liu, C.~Lageman, B.~Anderson, and G.~Shi, ``An arrow-hurwicz-uzawa type flow
  as least squares solver for network linear equations,'' \emph{arXiv preprint
  arXiv:1701.03908}, 2017.

\bibitem{Chen_2002}
H.-F. Chen, \emph{Stochastic approximation and its applications}.\hskip 1em
  plus 0.5em minus 0.4em\relax Springer Science \& Business Media, 2006,
  vol.~64.

\end{thebibliography}
\end{document}